\documentclass[11pt]{amsart}
\usepackage{amssymb}
\usepackage{amscd}
\usepackage{mathrsfs}
\usepackage{hyperref}

\newtheorem{mthm}{Theorem}

\newtheorem{thm}{Theorem}[section]
\newtheorem{cor}[thm]{Corollary}
 
\newtheorem{clm}[thm]{Claim}
\newtheorem{lem}[thm]{Lemma}

\theoremstyle{definition}
\newtheorem{defn}[thm]{Definition}
\newtheorem{defns}[thm]{Definitions}
\newtheorem{exmp}[thm]{Example}

\newtheorem{rem}[thm]{Remark}

\theoremstyle{remark}

\numberwithin{equation}{section}

\newcommand{\al}{\alpha}
\newcommand{\be}{\beta}
\newcommand{\ga}{\gamma}
\newcommand{\de}{\delta}
\newcommand{\pa}{\partial}   
\newcommand{\eps}{\epsilon}

\newcommand{\ph}{\phi}

\newcommand{\et}{\eta}

\newcommand{\la}{\lambda}
\newcommand{\rh}{\rho}
\newcommand{\si}{\sigma}
\newcommand{\ta}{\tau}

\newcommand{\ps}{\psi}

\newcommand{\GA}{\Gamma}
\newcommand{\LA}{\Lambda}
\newcommand{\DE}{\Delta}

\def\subs{\subseteq}
\def\normal{\triangleleft}
\def\iso{\cong}
\def\iff{{\quad\Longleftrightarrow\quad}}

\def\bs{{\backslash}}

\newcommand{\overto}[1]{{\buildrel{#1}\over\longrightarrow}}

\newcommand{\Set}[2]{{ \left\{ {#1} \,:\, {#2}\right\} }}

\newcommand{\ol}[1]{\overline{#1}}
\newcommand{\ben}{\begin{enumerate}}
\newcommand{\een}{\end{enumerate}}
\newcommand{\beit}{\begin{itemize}}
\newcommand{\enit}{\end{itemize}}
\newcommand{\bqr}{\begin{eqnarray}}
\newcommand{\eqr}{\end{eqnarray}}
\newcommand{\beq}{\begin{equation}}
\newcommand{\eeq}{\end{equation}}
\newcommand{\bqrn}{\begin{eqnarray*}}
\newcommand{\eqrn}{\end{eqnarray*}}
\newcommand{\lm}{\lim_{n\to\infty}}

\def\subs{\subseteq}
\def\GL{{\bf\rm GL}}

\def\SL{{\bf\rm SL}}
\def\PSL{{\bf\rm PSL}}
\def\SO{{\bf\rm SO}}

\def\SU{{\bf\rm SU}}

\def\supp{{\rm supp}}
\def\dist{{\rm dist}}

\def\locfld{{\mathbf F}}
\def\bH{{\mathbf H}}
\def\bR{{\mathbf R}}
\def\bC{{\mathbf C}}
\def\bQ{{\mathbf Q}}
\def\bZ{{\mathbf Z}}
\def\bN{{\mathbf N}}
\newcommand{\PrjSp}[2]{{ {{#1}}P^{{#2}} }}

\def\locfld{{k}}

\def\Id{{\it I\! d}}

\def\sB{{\mathcal{B}}}

\def\sF{{\mathcal{F}}}

\def\sP{{\mathcal{P}}}

\def\Aut{{\rm Aut}\,}

\def\Out{{\rm Out}\,}

\def\Isom{{\rm Isom}}
\def\Homeo{{\rm Homeo}}
\def\Hor{{\rm Hor}}
\def\hor{{\rm hor}}

\def\rmand{{\qquad{\rm and}\qquad}}

\def\rmwhere{{\qquad{\rm where}\qquad}}

\newcommand{\myemph}[1]{{\emph{#1}}}

\title[Measurable rigidity on homogeneous spaces]
{Measurable rigidity of actions on\\ infinite measure homogeneous spaces, II} 

\author{Alex Furman}\thanks{Supported in part by NSF grant DMS-0094245, 
and BSF USA-Israel grant 2004345.}

\address{Mathematics Statistics and Computer Science, University of Illinois at Chicago.}
\email{furman@math.uic.edu}


\begin{document}

\begin{abstract}
We consider the problems of measurable isomorphisms and joinings,
measurable centralizers and quotients for certain classes
of ergodic group actions on \emph{infinite measure} spaces.
Our main focus is on systems of algebraic origin: actions of lattices and 
other discrete subgroups $\Gamma<G$ on homogeneous spaces 
$G/H$ where $H$ is a sufficiently rich unimodular subgroup 
in a semi-simple group $G$.
We also consider actions of discrete groups of isometries $\Gamma<\Isom(X)$ 
of a pinched negative curvature space $X$, acting on the space of horospheres 
$\Hor(X)$.
For such systems we prove that the only measurable isomorphisms, joinings, 
quotients etc. are the obvious algebraic (or geometric) ones. 
This work was inspired by the previous work of Shalom and Steger,
but uses completely different techniques which lead to more general results.
\end{abstract}

\maketitle


\section{Introduction and Statement of the Main Results}
 
The starting point of our discussion is the following beautiful result of 
Yehuda Shalom and Tim Steger: 
\begin{thm}[Shalom-Steger, \cite{ShalomSteger}]\label{T:SL2}
Let $\Gamma_{1}, \Gamma_{2}$ be two abstractly isomorphic lattices
in $\SL_{2}(\bR)$, and $\tau:\Gamma_{1}\overto{\iso}\Gamma_{2}$ be 
the isomorphism.
Then the only possible measurable isomorphisms between the linear actions  
of $\Gamma_{i}$ on $\bR^{2}$ are algebraic ones.
More precisely, if $T:\bR^{2}\to\bR^{2}$ is a strinctly measure class preserving map 
with $T(\gamma x)=\gamma^{\ta} T(x)$  a.e. $x\in\bR^{2}$ and all $\gamma\in\Gamma_{1}$, then there exists $A\in \GL_{2}(\bR)$ so that 
$\gamma^{\ta}=A\,\gamma A^{-1}$ for all $\gamma\in\Gamma_{1}$ and 
$T(x)=Ax$ a.e. on $\bR^{2}$.
\end{thm}
Obsereve that the linear $\SL_{2}(\bR)$-action on $\bR^{2}-\{0\}$ (which is measure-theoretically 
indistinguishable from $\bR^{2}$) is precisely
the action of $G=\SL_{2}(\bR)$ by left translations on the homogeneous space
$G/H$ where $H$ is the \emph{horocyclic} subgroup
\[
	H=\Set{\begin{pmatrix} 1 & s\\ 0 & 1\end{pmatrix}}{s\in\bR}.
\]
The action of $\Gamma$ on $G/H$ is closely related to its ``dual`` dynamical system 
-- the action of $H$ on $G/\Gamma$.
The latter is an algebraic description of the horocycle flow on the unit tangent bundle 
$S M$ to the Riemann surface $M=\bH^{2}/\Gamma=K\bs G/\Gamma$
(here we assume that $\Gamma$ is torsion free).
In the 1980s Marina Ratner discovered remarkable measurable rigidity properties of 
the horocycle flow, proving that all measurable isomorphisms \cite{Ratner:isom-hor:82},
measurable quotients \cite{Ratner:Factors:82}, and finally all joinings 
\cite{Ratner:Joinings:83} 
of such flows are algebraic, i.e. $G$-equivariant.
Shalom-Steger result above can be viewed as a ``dual companion'' of Ratner's
isomorphism theorem  \cite{Ratner:isom-hor:82}.
However, it is important to emphasize that despite the formal similarities theorem~\ref{T:SL2} 
does not seem to be directly related to (neither to imply, nor to follow from) any of the above results
of Ratner; it also cannot be deduced from the celebrated Ratner's Classification of invariant measures theorem \cite{Ratner:ICMtalk:94}, which contains \cite{Ratner:isom-hor:82}, \cite{Ratner:Factors:82}, 
\cite{Ratner:Joinings:83} as particular cases.

Shalom and Steger prove their theorem~\ref{T:SL2} (and other rigidity results, 
such as \ref{C:SLn} (1), (2) below)
ingeniously using unitary representations techniques. 
The present paper grew out of an attempt to give an alternative, purely dynamical
proof for this theorem and other related results from \cite{ShalomSteger}
\footnote{Hence the numeral II in the title of this paper.}.
The technique that has been developed for this purpose --- the alignment property, 
is quite a general and turns out to be very powerful.
We develop it as an abstract tool and show how to apply it to 
homogeneous spaces and spaces of horospheres. 
In the present work the rigidity phenomena of Theorem~\ref{T:SL2} are generalized
in several directions: (1) we consider homogeneous spaces of all semi-simple groups,
(2) we also consider spaces of horospheres in variable pinched negative curvature, 
(3) in the context of negative curvature we treat actions of discrete groups $\Gamma$ 
which are not necessarily lattices; 
(4) in all the above examples we prove rigidity results not only for isomorphisms, but
also for quotients and more generally for joinings.

\bigskip

Before stating the results we need to set a few conventions concerning ${\rm II}_{\infty}$ 
\emph{actions} -- these are measure preserving, ergodic group actions on non-atomic \emph{infinite measure} Lebesgue spaces (see section \ref{S:prelim} for more details). 
\begin{defns} \label{D:main}
Let $(X_{i},m_{i},\Gamma)$, $i=1,2$ be two ergodic measure-preserving actions
of a fixed group $\Gamma$ on infinite measure Lebesgue spaces $(X_{i}, m_{i})$,
and $\tau\in\Aut\Gamma$ be a group automorphism. 
A \myemph{morphism} or ($\tau$-twisted) \myemph{quotient map} between such systems 
is a measurable map $\pi:X_{1}\to X_{2}$ such that 
\[
	\pi_{*}m_{1}=const\cdot m_{2}\qquad\text{and}\qquad 
	\pi(\gamma x)=\gamma^{\tau} \pi(x)
\]
for all $\gamma\in\Gamma$ and $m_{1}$ a.e. $x\in X_{1}$. 
In particular, the first condition implies that the preimage $\pi^{-1}(E)$ of a set 
$E\subset X_{2}$ of finite $m_{2}$-measure has finite $m_{1}$-measure.
A ($\tau$-twisted) \myemph{isomorphism} is a measurable bijection with the same properties.
An (${\rm Id}\times\tau$-twisted) \myemph{joining} of such 
systems is a measure $\bar{m}$ on 
$X=X_{1}\times X_{2}$ which is invariant under the (twisted) diagonal $\Gamma$-action 
\[
	\gamma:(x_{1},x_{2})\mapsto (\gamma x_{1},\gamma^{\tau}x_{2})
\]
and such that the projections $\pi_{i}:X\to X_{i}$, $i=1,2$ are morphisms, i.e. 
\[
	(\pi_{i})_{*} \bar{m}=const_{i}\cdot m_{i}\qquad(i=1,2).
\]
Two systems are \myemph{disjoint} if they admit no joinings\footnote{Note that for infinite measure systems the product measure $m_{1}\times m_{2}$ is a not a joining}.
Given an infinite measure preserving ergodic system $(X,m,\Gamma)$
its measurable \myemph{centralizer} is defined to be the group of all measurable
(possibly twisted) automorphisms of the $\Gamma$-action on $(X,m)$. 
Similarly \myemph{self joinings} are joinings of $(X,m,\Gamma)$ with itself. 
\end{defns}

This framework of ${\rm II}_{\infty}$ actions in many respects parallels that of 
${\rm II}_{1}$ actions -- the classical theory of ergodic \emph{probability measure} 
preserving actions.
For example self joinings control centralizers and quotients of a given system, 
and joinings between two systems control isomorphisms and common quotients
(see section~\ref{S:prelim} for details). 
Hence stating the results below we consider separately: 
\begin{itemize}
\item[(I)] results about a single action for which centralizers, quotients and self-joinings are studied;
\item[(II)] results on pairs of systems, for which isomorphisms, common quotients and
general joinings are analyzed.  
\end{itemize}


We shall be mostly interested in actions of discrete subgroups $\GA< G$ on homogeneous
spaces $X=G/H$, where $G$ is a locally compact (always second countable) group, 
$H$ is a closed subgroup so that $G/H$ carries an infinite $G$-invariant measure 
$m_{G/H}$. 
The main results of the paper assert that under certain mild conditions:
\begin{itemize}
\item[(I)] 
measurable centralizers, quotients and self joinings of a single system $(G/H,\Gamma)$, 
\item[(II)] 
measurable isomorphisms and joinings between two such systems $(G_{i}/H_{i},\Gamma)$
$i=1,2$
\end{itemize}
are \emph{algebraic}, i.e., essentially coincide with the corresponding notions (centralizers, quotients, isomorphisms, joinings, etc.) for the transitive 
$G$-action on $X=G/H$.
Let us describe algebraic centralizers and quotients more explicitly:
\begin{exmp}[Algebraic Centralizers and Quotients] \label{E:algebraicCJQ}
Let $X=G/H$ be a homogeneous space with an infinite $G$-invariant measure $m$.
For the transitive $G$-action on $(X,m)$ we have the following:
\begin{description}
\item[Centralizers]
    The centralizer of the $G$-action on $X$ in both measurable and set-theoretic
    sense is the group $\LA=N_G(H)/H$, 
    where $\lambda=n_{\lambda}H\in\Lambda$ acts on 
    $X=G/H$ by $\lambda: x=gH\mapsto \lambda x=gn_{\lambda}H$. 
\item[Quotients]
   Any $G$-equivariant measurable quotient of $G/H$ is $G/H^{\prime}$ via
   $\pi:gH\mapsto gH^{\prime}$, where $H<H^{\prime}$ and $H^{\prime}/H$ 
   is compact.
\end{description}
\end{exmp}

\bigskip

We shall consider semi-simple Lie groups $G$ and a class of 
closed unimodular subgroups $H<G$ which we call ``suitable''
(see definition~\ref{D:suitableH}). For rank one real Lie group $G$
suitable subgroups are any closed subgroup $N<H<MN$ where $N$ 
is the horospherical subgroup and $M<K$ is the centralizer of the Cartan $A$ in $K$. 
The assumptions on $\Gamma<G$ will vary: requiring $\Gamma<G$ to be 
a lattice would be sufficient to establish rigidity for all the examples;  
for homogeneous spaces $G/H$ of rank one real Lie groups $G$, 
a wider class of discrete subgroups $\Gamma<G$ can be shown to be rigid. 
We start with these latter cases.

\subsection*{Homogeneous spaces of rank one real Lie groups}
\begin{mthm}[{\bf Real Rank One: Centralizers, Self Joinings and Quotients}]
\label{T:self-rk1}
Let $G$ be a real, connected, Lie group of rank one with trivial center, 
$N< G$ its horospherical subgroup,  and $H< G$ a proper closed 
unimodular subgroup $N<H<MN$.
Suppose that $\GA< G$ is a discrete subgroup acting ergodically on 
$(G/N,m_{G/N})$ and hence also on the homogeneous space $(X,m)=(G/H,m_{G/H})$. 

Then the $\GA$-action on $(X,m)$ has only algebraic 
centralizers and quotients as described in Example~\ref{E:algebraicCJQ}, 
and any ergodic self-joining descends to an algebraic centralizer of the
algebraic quotient $G/MN$.
\end{mthm}
Let us describe the scope of this theorem.
The possibilities for $H<G$ as in the theorem are quite restricted:  
the spaces $G/H$ are compact extensions of $G/MN$ -- the space
of horospheres $\Hor(S\bH)$ in the unit tangent bundle $S\bH$ to the symmetric space
$\bH\iso G/K$ of $G$. 
However, the condition on a discrete subgroup $\Gamma< G$ is quite mild. 
Examples of such subgroups include
\begin{itemize}
\item
	Any lattice $\GA$ in $G$ (both uniform and non-uniform ones). 
	Ergodicity of the $\Gamma$ action on $(G/H,m_{G/H})$ follows from 
	Moore's ergodicity theorem ($H$ is not precompact and hence acts ergodically
	on $G/\Gamma$).
\item
	Let $\LA< G$ be a lattice and $\GA\normal \LA$ so that $\LA/\GA$ is nilpotent.
	Then $\GA$ acts ergodically on $G/MN=\Hor(S\bH)$.
	This was proved by Babillot and Ledrappier \cite{BabillotLedrappier:TIFR:98} 
	for the case where 
	$\LA/\GA$ is Abelian and 	$\LA< G$ is a uniform torsion free lattice.
	In \cite{Kaimanovich:SAT:02} Kaimanovich showed that in this context ergodicity 
	of the $\GA$-action on the space of horospheres 
	is equivalent to the ergodicity of the $\GA$-action on the sphere at infinity 
	$\pa\bH=G/MAN$ which, in turn, is equivalent to the lack of non-constant 
	bounded harmonic functions on the regular cover $\bar{M}=\GA\bs\bH$
	of the finite volume manifold $M=\LA\bs\bH$. 
	For nilpotent covering group $\LA/\GA$ the latter is well known 
	(e.g. Kaimanovich \cite{Kaimanovich:erghor:00}).
\end{itemize}

\medskip

The next rigidity result requires a stronger assumption on a discrete group $\GA$ 
in a rank one Lie group $G$.
Let $\bH=G/K$ denote the symmetric space of $G$, and $\pa \bH=G/MAN$ its boundary. 
\begin{defn}\label{D:E2}
We shall say that a torsion free $\GA$ satisfies condition (E2) if the following equivalent conditions hold
(the equivalence is due to Sullivan \cite{Sullivan:BAMS:82}) 
\begin{itemize}
\item[{\rm (E2a)}]
	$\GA$ acts ergodically on $\pa\bH\times\pa\bH$ with respect to 
	the standard measure class (that of the $K\times K$-invariant measure).
\item[{\rm (E2b)}] 
	The geodesic flow is ergodic on $S\bH/\GA$.
\item[{\rm (E2c)}]
	The Poincar\'e series $\sum_{\ga \in\GA} e^{-s\cdot d(\ga p,p)}$ diverges at 
	$s=\de(\bH)$ where 
	\[
		\de(\bH)=\lim_{{R\to\infty}} \frac{1}{R}\log{\rm Vol}(B(p,R))
	\]
	denotes the volume growth rate of the symmetric space $\bH$.  
\end{itemize}
\end{defn}
These conditions are satisfied by any lattice $\GA< G$. 
In \cite{Guivarch:ETDS:89} Guivarc'h considered geodesic flows on regular covers 
of compact hyperbolic manifolds. 
His results (extending previous work of M. Rees) in particular imply that if $\Lambda<G$
is a uniform lattice and $\GA\normal\LA$ then $\Gamma$ satisfies (E2b) iff  
a simple random walk on $\LA/\GA$ is recurrent, which occurs iff $\LA/\GA$
is a finite extension of $\bZ^{d}$ with $d\le 2$.

\begin{mthm}[{\bf Real Rank One: Rigidity of actions}]
\label{T:rigidity-rk1}
Let $G_1, G_2$ be real, connected, non compact,  rank one Lie groups with trivial centers, 
$N_i< G_i$ the horospheric subgroups, $H_i< G_i$ closed unimodular
subgroups with $\check{H}_{i}=N_{i}<H_{i}<\hat{H}_{i}=M_{i}N_{i}$, 
and $(X_i,m_i)=(G_i/H_i,m_{G_i/H_i})$. 
Let $\GA_i< G_i$ be discrete subgroups satisfying condition {\rm (E2)}
and acting ergodically on $(G_{i}/N_{i},m_{G_{i}/N_{i}})$.
Assume that $\GA_1$ and $\GA_2$ are isomorphic as abstract groups $\ta:\GA_1\overto{\iso}\GA_2$ 
and that $(X_i,m_i,\GA_i)$ admit a $\ta$-twisted (ergodic) joining $\bar{m}$. 
Then 
\begin{enumerate}
\item
	Then $\tau:\Gamma_{1}\to\Gamma_{2}$ extends to an isomorphism of 
	the ambient groups $\tau:G_{1}\overto{\iso}G_{2}$
	which maps $N_{1}$ onto $N_{2}$.
\item
	The joining $\bar{m}$ descends to an algebraic isomorphism 
	$T^{\prime}:G_{1}/\hat{H}_{1}\to G_{2}/\hat{H}_{2}$
	between algebraic quotients of the original systems,
	with $\hat{H}_{i}=M_{i}N_{i}> H_{i}$.
\item
	If the actions not only admit an ergodic joining, but actually a measurable $\ta$-twisted 
	isomorphism $T:X_{1}\to X_{2}$, then the isomorphism $\ta:G_{1}\to G_{2}$
	as in {\rm (1)} in addition maps $H_{1}$ onto $H_{2}$ and for some 
	$\la\in N_{G_{2}}(H_{2})$
	we have almost everywhere
	$T(gH_{1})=\la g^{\ta} H_{2}.$
\end{enumerate}
\end{mthm}
The above Theorem in particular applies to $\GA_i< G_i$ being lattices.
However due to Mostow Rigidity the only examples of abstractly isomorphic but not conjugate 
lattices occur in $G_1=G_2=\PSL_2(\bR)$. 
In these cases $X_1=X_2=(\bR^2-\{0\})/x\sim \pm x$, but an easy
modification gives a similar rigidity result for $\Gamma_{i}<\SL_{2}(\bR)$ acting linearly on $\bR^{2}$.
In this case statement (3) gives Shalom-Steger's result~\ref{T:SL2}.

In addition to lattices in $\SL_{2}(\bR)$ there are many examples of infinite 
covolume discrete subgroups in rank one $G$ satisfying condition (E2), in particular
these examples include certain normal subgroups in uniform lattices, namely 
fundamental groups of $\bZ$ or $\bZ^{2}$ regular covers
of compact locally symmetric spaces.
This opens the possibility for the same group $\Gamma$ to be embedded
as a discrete subgroup satisfying (E2) in different rank one groups $G_{1}$ and $G_{2}$.
It is known that most (conjecturally all) arithmetic lattices in rank one groups $G\simeq \SO_{n,1}$ 
and $\SU_{n,1}$ have a finite  index subgroup $\Lambda$ with infinite abelianization
(cf. \ref{Lubotzky:BettiTau:96}, \ref{RaghunathanVenkataramana:Betti:92}). 
In particular it would fit in an exact sequence $\Gamma\to\Lambda\to\bZ$, 
and often the kernel $\Gamma$ is expected to be a free group on infinitely many generators 
$F_{\infty}$.

Another class of examples is obtained by embedding the fundamental group 
$\GA=\pi_1(S)$  of a closed orientable surface $S$ of genus $g\ge 2$ in $\PSL_{2}(\bC)$.
Let $\ph$ be a pseudo-Anosov diffeomorphism of $S$, and $\GA_\ph=\bZ\ltimes_{[\ph]}\GA$ 
denote the semi-direct product defined by $[\ph]\in\Out\GA$. Then $\GA_\ph$ is a fundamental
group for the mapping torus $3$-manifold $M_\ph=S\times [0,1]/(x,0)\sim(\ph(x),1)$.
By the famous hyperbolization theorem of Thurston (\cite{Otal:hyperbolization:01}) 
such an $M_\ph$ admits a 
hyperbolic structure, i.e. $\GA_\ph$ is a cocompact lattice in $G=\Isom_+(\bH^3)=\PSL_2(\bC)$.
It contains the surface group $\GA$ as a normal subgroup with $\GA_\ph/\GA\iso\bZ$.
Thus $\GA$ satisfies (E2).
Therefore surface groups $\GA$ can appear as a discrete subgroup with condition (E2) 
in a variety of ways in $\PSL_2(\bR)\iso\Isom_{+}(\bH^{2})$ and in $\PSL_{2}(\bC)\iso\Isom_{+}(\bH^{3})$.
In the former there is a continuum of such embeddings -- parametrized by the Teichmuller space; 
in the latter there are (at least) countably many such embeddings 
defined by varying a pseudo-Anosov element $\ph$ in the Mapping Class Group of $S$.

\begin{rem}
The rank one results can be extended to the geometric context of the spaces of horospheres
in manifolds of variable negative curvature.
We formulate this result (Theorem~\ref{T:rigidity-hor}) in Section~\ref{S:Horospheres}.
\end{rem}

\medskip

\subsection*{Homogeneous spaces of general semi-simple groups}

It turns out that rigidity phenomena for actions of \emph{lattices} are quite wide spread among 
homogeneous spaces $G/H$ of semisimple groups $G$ and sufficiently large unimodular $H<G$. 
Before formulating general results consider the following: 
\begin{exmp}\label{E:SLn}
Let $k$ be a local field, e.a. $\bR$, $\bC$, or a finite extension of $\bQ_p$ for a prime $p$, 
and let $G=\SL_n(\locfld)$. 
Then $X=\locfld^n\setminus\{0\}$ is the homogeneous space 
$G/H$ for the unimodular closed subgroup
\[
	H=\Set{ g\in\SL_n(\locfld) }{g_{11}=1,\ g_{21}=\dots=g_{n1}=0}.
\]
More generally, given  a partition $n=n_1+\cdots+n_m$ (with $m>1$ and $n_{i}\in\bN$) 
consider the subgroup $Q<G=\SL_n(\locfld)$ consisting of the upper triangular block matrices of the form 
\begin{equation}\label{e:blocks}
	\begin{pmatrix}
         A_{11} & B_{12} & \cdots & B_{1m}   \\
	0 & A_{22} & \cdots & B_{2m} \\
       	\vdots & \vdots & \ddots& \vdots\\
	0 & \cdots  & 0 & A_{mm}
\end{pmatrix}
\end{equation}
where $A_{ii}\in\GL_{n_i}(\locfld)$ and $B_{ij}\in M_{n_i\times n_j}(\locfld)$. 
Let $\check{H}\normal Q$ denote the closed subgroup consisting of matrices with
$\det A_{ii}=1$, ($i=1,\dots,m$); 
and let $\hat{H}\normal Q$ denote a slightly larger subgroup 
consisting of block matrices as above with $|\det A_{ii}|=1$, 
where $|\cdot|:\locfld\to[0,\infty)$ denotes the absolute value on $\locfld$.
Then $\check{H}<\hat{H}$ and any intermediate closed subgroup 
$\check{H}<H<\hat{H}$ are unimodular and $X=G/H$ carries an infinite 
$G$-invariant measure $m=m_{G/H}$. 
\end{exmp}
We shall now describe the most general setting for our rigidity results:
$G$ will be a semi-simple group (in a very general sense), while $H<G$ will be restricted
to some rather class of subgroups, which for a lack of a better term we call ``\emph{suitable}''. 
Subgroups $H<G$ which appear in Theorems~\ref{T:self-rk1} and \ref{T:rigidity-rk1}
and in Example~\ref{E:SLn} above are suitable. 
The formal definition/construction is the following: 
\begin{defn}[Suitable subgroups $H$ is semi-simple $G$]
\label{D:suitableH}
Let $A$ be a finite set, for $\alpha\in A$ let $k_{\alpha}$ be a local field of zero characteristic
and ${\mathbf G}_{\alpha}$ be some connected, semi-simple linear algebraic $\locfld_{\alpha}$-group. 
The product  
\begin{equation}
\label{e;defG}
	G=\prod_{\alpha\in A}{\mathbf G}_{\alpha}(\locfld_{\alpha})
\end{equation}
of $\locfld_{\alpha}$-points of the corresponding $\locfld_{\alpha}$-groups taken with the Hausdorff topology is a localy compact second countable group. 
We shall refer to such groups as just \emph{semi-simple}. 
Groups $H<G$ appearing in the following construction will be called \emph{suitable}. 

For each $\alpha\in A$ choose a parabolic subgroup 
${\mathbf Q}_{\alpha}<{\mathbf G}_{\alpha}$
and, taking the product of $\locfld_{\alpha}$-points of these groups, 
form a closed subgroup  $Q= \prod_{\alpha\in A} {\mathbf G}_{\alpha}(\locfld_{\alpha})$ 
in $G$.
To any such subgroup  $Q<G$, we shall them  \emph{parabolic}, 
we associate two closed unimodular  
subgroups $\check{H}\normal\hat{H}<G$ with $\check{H}<[Q,Q]<\hat{H}$ and
$\hat{H}/\check{H}$ compact. 
Any intermediate closed subgroup $\check{H}<H<\hat{H}$ 
will be called \emph{suitable}.
The subgroups $\check{H}<\hat{H}$ are constructed as follows:
let ${\mathbf Q}_{\alpha}={\mathbf R}_{\alpha}\ltimes {\mathbf V}_{\alpha}$ be the Levi decomposition
into the the unipotent radical and a reductive factor, and let ${\mathbf S}_{\alpha}$ denote the product 
of all the $\locfld_{\alpha}$-isotropic factors of the the commutator subgroup 
$[{\mathbf R}_{\alpha}, {\mathbf R}_{\alpha}]$. 
Hence  ${\mathbf S}_{\alpha}$ is a semi-simple $\locfld_{\alpha}$-group, and we let 
$\check{H}=\prod_{\alpha\in A} ({\mathbf S}_{\alpha}\ltimes {\mathbf V}_{\alpha})(\locfld_{\alpha})$.
It is a closed cocompact subgroup of the commutator 
$
	[Q,Q]=\prod_{\alpha\in A} ([{\mathbf R}_{\alpha},{\mathbf R}_{\alpha}]
		\ltimes {\mathbf V}_{\alpha})(\locfld_{\alpha}).
$
The abelianization $Q/[Q,Q]$ of $Q$ is the product of $\locfld_{\alpha}$-torii 
$\prod \GL_{1}(\locfld_{\alpha})^{r_{\alpha}}$, 
(here $r_{\alpha}$ is the $\locfld_{\alpha}$-rank of ${\mathbf R}_{\alpha}$). 
These torii have unique maximal compact subgroups 
$U_{\alpha}=\Set{ (t_{1},\dots,t_{r_{\alpha}}) }{ |t_{1}|_{\alpha}=\dots= |t_{r_{\alpha}}|_{\alpha}=1 }$. 
We denote by $\hat{H}$ the preimage of the product $\prod_{\alpha\in A} U_{\alpha}$ under the Abelianization epimorphism $Q\to Q/[Q,Q]$.
\end{defn}
\begin{rem}
\label{R:suitable}
In a given semi-simple group $G$, the collection of all suitable subgroups $H<G$ is
divided into families of groups related to a given parabolic $Q<G$, in each such family
the groups $H$ share a common cocompact normal subgroup $\check{H}$ and 
a common compact extension $\hat{H}$ with $[Q,Q]$.
All ``suitable'' subgroup $H<G$ are closed unimodular and  contain the maximal
unipotent subgroup $N<G$ ($N$ is the unipotent radical of the minimal parabolic $P<G$);
it seems that the converse is also true: any unimodular closed subgroup containing $N$ 
appears in the construction ~\ref{D:suitableH}.
\end{rem}
  
\medskip

\begin{mthm}[{\bf General case: Centralizers, Self Joinings and Quotients}] 
\label{T:self-higher}
Let $G=\prod {\mathbf G}_{\alpha}(\locfld_{\alpha})$ be a semi-simple group and $H<G$ 
be a suitable subgroup as in \ref{D:suitableH} associated to 
a parabolic $Q<G$. Let $\Gamma<G$ be a lattice acting by left translations on
$(X,m)=(G/H,m_{G/H})$.

Then the only measurable centralizers and quotients of the $\Gamma$-action on $(X,m)$
are algebraic (as in \ref{E:algebraicCJQ}), and any ergodic self-joining descends to an
algebraic automorphism of an algebraic quotient 
$(\hat{X},\hat{m})=(G/\hat{H},m_{G/\hat{H}})$, and is itself a quotient of an algebraic
automorphism of an algebraic extension 
$(\check{X},\check{m})=(G/\check{H},m_{G/\check{H}})$.
\end{mthm}
\begin{cor}\label{C:SLn}
Let $\locfld$ be a local field, $\GA< \SL_n(\locfld)$ be a lattice and let $(X,m)$ denote
the vector space $\locfld^n$ with the Lebesgue measure, with the linear $\GA$-action. 
Then 
\begin{enumerate}
\item 
	the measurable centralizer of the system $(X,m,\GA)$ 
	consists only of homotheties: $\ x\mapsto \la x$, where $\la\in\locfld^{*}$;
\item
	the only measurably proper quotients of $(X,m,\GA)$ are of the
	form $\locfld^n/C$ where $C$ is a closed subgroup of the compact 
	Abelian group of $\locfld$-units $U_{\locfld}=\{u\in\locfld\,:\, |u|=1\}$;
\item
	the only ergodic self joinings are on graphs of homotheties 
	\[
		\{(x,\la x)\mid x\in\locfld^n\}\qquad (\lambda \in \locfld^{*}).
	\]
\end{enumerate}
\end{cor}
\begin{rem}
Items (1) and (2) in the above Corollary were first proved by Shalom-Steger \cite{ShalomSteger}.
\end{rem}

For higher rank groups our techniques are restricted to lattices.
Due to Mostow-Margulis rigidity we are not able to vary the embedding
of a given lattice in a higher rank group $G$. 
Hence we shall consider a fixed lattice $\Gamma< G$,
but will still be able to vary the homogeneous space $G/H$.
\begin{mthm}[{\bf General case: Rigidity for actions}] \label{T:rigidity-higher}
Let $G=\prod {\bf G}_i(\locfld_i)$ be a semi-simple group, $\Gamma<G$ a lattice, 
and $H_1, H_2< G$ be two suitable subgroups as in \ref{D:suitableH}.
Assume that the $\GA$-actions on $(X_{i}, m_{i})=(G/H_{i},m_{G/H_{i}})$ admit
an ergodic joining $\bar{m}$. 

Then $X_{1}$, $X_{2}$ share a common algebraic quotient 
$(\hat{X},\hat{m})=(G/\hat{H},m_{G/\hat{H}})$, and are common algebraic extension
$(\check{X},\check{m})=(G/\check{H},m_{G/\check{H}})$; 
the joining $\bar{m}$ descends to an algebraic automorphism of the
$\Gamma$-action on $(\hat{X},\hat{m})$ and is a quotient of an algebraic
automorphism of $(\check{X},\check{m})$. 

Furthermore, if the original $\Gamma$-actions on $(X_{i},m_{i})$ are isomorphic, 
say via $T:X_{1}\to X_{2}$, then for some $q\in Q$ 
\[
	H_{1}= qH_{2}q^{-1}\qquad\text{and}\qquad T(gH_{1})=gqH_{2}
\]
for $m_{1}$-a.e. $gH_{1}\in X_{1}$.
\end{mthm}

\medskip

\subsection*{Organization of the paper}
In section~\ref{S:prelim} we summarize some basic properties of our setup, 
such as joinings of infinite measure preserving group actions.
In section \ref{S:generalCJQ} we introduce the notion of
\emph{alignment systems} and analyze centralizers, self-joinings and quotients
of general actions admitting alignment systems. 
This general theory continues in Section~\ref{S:boundary} where joinings
of two general systems are studied under the assumption that each has 
an alignment property.
In section \ref{S:alignment} we prove the alignment property for our main examples:
homogeneous spaces and spaces of horospheres. 
In section \ref{S:Horospheres} joinings between spaces of horospheres are 
studied. In section~\ref{S:proofs-rigidity} the final results on homogeneous 
spaces are proved.

\subsection*{Acknowledgments}
This paper has been inspired by the work of Yehuda Shalom and Tim Steger \cite{ShalomSteger}.
I would like to thank them for sharing their results and insights.
I would also like to that Gregory Margulis and Shahar Mozes for stimulating discussions 
about these topics, and Marc Burger and FIM at ETH Zurich for their hospitality
during summer 2002.

\section{Preliminaries} \label{S:prelim}

\subsection{Strictly Measure Class Preserving maps}

We first discuss some technical points regarding ${\rm II}_{\infty}$ actions. 
Let $(X_i,m_i,\GA)$, $i=1,2$  be two such actions of a fixed group $\Gamma$.
A measurable map $T:(X_1,m_{1})\to (X_2,m_{2})$ will be called 
\myemph{strictly measure class preserving} if $\ T_{*}m_{1}\sim m_{2}\ $ and the Radon-Nikodym
derivative $T^{\prime}$ is almost everywhere positive and finite. 
Note that such a map as the projection $\bR^{2}\to \bR$, $(x,y)\mapsto x$, is not  \emph{strictly}
measure class preserving, although is usually considered measure class preserving.

Any strictly measure class preserving $\GA$-equivariant map $X_{1}\to X_{2}$ 
between ergodic measure preserving $\GA$-actions, has a $\GA$-invariant, 
and hence a.e. constant positive and finite Radon-Nikodym derivative.
Therefore
\begin{lem}
If $(X_i,m_i,\GA)$, $i=1,2$ are two ergodic infinite measure preserving systems, 
and some fixed $\ta\in\Aut\Gamma$, then:
\begin{enumerate}
\item
	Any strictly measure class preserving map $T:X_{1}\to X_{2}$ 	satisfying 
	\[
		T(\gamma x)=\gamma^{\ta}T(x)\qquad (\gamma\in \Gamma)
	\]
	for a.e. $x\in X_{1}$, 	is a ($\ta$-twisted) quotient map.
	If furthermore $T$ is invertible then $T$ is a ($\tau$-twisted) isomorphism.
\item
	A measure $\bar{m}$ on $X=X_{1}\times X_{2}$ invariant under 
	the (${\rm Id}\times\tau$-twisted) diagonal $\Gamma$-action
	\[
		\ga:(x_1,x_2)\mapsto (\ga x_1,\ga^\ta x_2)
	\]
	for which the projections $\pi_{i}:(X,\bar{m})\to (X_{i},m_{i})$ are strictly 
	measure class preserving,
	is a ($\ta$-twisted) joining of $(X_{i}, m_{i},\Gamma)$, $i=1,2$.
\end{enumerate}
\end{lem}

\medskip

\begin{lem}
\label{L:quotient-by-comp}
Let $(X,m,\Gamma)$ be an ergodic, infinite measure preserving system, and 
suppose that $L$ is a locally compact group with a faithful measurable action 
by measure class preserving transformations on $(X,m)$, commuting with $\Gamma$. 
Then
\begin{enumerate}
\item
	There exists a continuous multiplicative character $\Delta:L\to \bR_{+}^{*}$ so that
	$g_{*}m=\Delta(g)m$ for all $g\in L$;
\item
	Any compact subgroup $K<L$ acts by measure preserving transformations on $(X,m)$,
	the $K$-orbits define a quotient system 
	\[
		P:(X,m,\Gamma)\to (X^{\prime}, m^{\prime}, \Gamma)\qquad
		\text{where}\qquad X^{\prime}=X/K,\ m^{\prime}=m/K
	\] 
	with $P:x\in X\mapsto Kx\in X^{\prime}$.
\end{enumerate}
\end{lem}
\begin{proof}
(1). The derivative cocycle  $\Delta(g,x)=\frac{dg_{*}m}{dm}(x)$ is a measurable map 
$L\times X\to \bR_{+}^{*}$ which is invariant under the $\Gamma$-action on $X$. 
By ergodicity, $\Delta(g,x)$ is a.e. constant $\Delta(g)$ on $X$.
Hence $\Delta$ is a measurable character $L\to \bR_{+}^{*}$.
It is well known that any measurable homomorphism between locally compact second countable groups
are continuous.

(2). The multiplicative positive reals $\bR_{+}^{*}$ do not have any non-trivial 
compact subgroups.
Thus $\Delta$ is trivial on $K$, i.e., $K$ acts by measure preserving transformations on $(X,m)$. 
The space of $K$-orbits  $X^{\prime}$ inherits: (1) 
a measurable structure from $X$ (because $K$ is compact), 
(2) the action of $\Gamma$  (because it commutes with $K$), and 
(3) the measure $m^{\prime}$, as required.
\end{proof}

\subsection{Joinings of ${\rm II}_{\infty}$ systems}\label{SS:Joinings}
Let us point out some facts about joinings of ergodic infinite measure preserving systems:
\begin{enumerate}
\item
   Any joining between two ergodic infinite measure preserving systems disintegrates into 
   an integral over a probability measure space of a family of ergodic joinings.
\item
   If $(X\times Y,\bar{m})$ is a $\ta$-twisted joining of ergodic infinite measure preserving systems
   $(X,m,\GA)$ and $(Y,n,\GA)$ (where $\bar{m}$ projects as $c_1\cdot m$ on $X$ and 
   as $c_2\cdot n$ on $Y$) then there exist unique up to null sets measurable maps
   $X\to \sP(Y)$, $x\mapsto \mu_x$, and $Y\to \sP(X)$, $y\mapsto \nu_y$, so that
   \[
   	\bar{m}=c_1\cdot \int_X \de_x\otimes\mu_x\,dm(x)=c_2\cdot \int_Y \nu_y\otimes\de_y\,dn(y)
   \] 
   Furthermore $\mu_{\ga x}=\ga_* \mu_x$ and $\nu_{\ga^\ta y}=\ga_*\nu_y$ for all $\ga\in\GA$ and
   $m$-a.e. $x\in X$ and $n$-a.e. $y\in Y$. 
\item
   In contrast to actions on probability spaces, group actions on infinite measure systems 
   do not always admit a joining\footnote{Note that the product measure 
   $m_1\times m_2$ is not a joining.}.
   Existence of a joining is an \emph{equivalence relation} between  
   ${\rm II}_{\infty}$ actions of a fixed group $\Gamma$. 
   Indeed, reflexivity and symmetry are rather obvious (as self joining use the diagonal
   $\int\delta_{x}\otimes\delta_{x}\,dm(x)$ measure). 
   For transitivity, given a joining $\bar{m}_{12}$ of $(X_{1},m_{1},\Gamma)$ with 
   $(X_{2},m_{2},\Gamma)$, 
   and $\bar{m}_{23}$ of $(X_{2},m_{2},\Gamma)$ with $(X_{3},m_{3},\Gamma)$
   one can use the decompositions
   \[
   	\quad
	\bar{m}_{12}=c_{12}\cdot \int_{X_{2}}\mu^{(1)}_{y}\otimes\delta_{y}\,dm_{2}(y),
	\quad
   	\bar{m}_{12}=c_{23}\cdot \int_{X_{2}}\delta_{y}\otimes\mu^{(3)}_{y}\,dm_{2}(y)
   \]
   with $\mu^{(i)}:X_{2}\to\sP(X_{i})$ ($i=1,3$) measurable functions, 
   in order to construct  the ``amalgamated'' joining $\bar{m}$ of $X_{1}$ with $X_{3}$
   by setting
   \[
   	\bar{m}=\int_{X_{2}}\mu^{(1)}_{y} \otimes \mu^{(3)}_{y}\,dm_{2}(y).
   \]
   Note that this amalgamated joining need not be ergodic, even if $\bar{m}_{12}$ and 
   $\bar{m}_{23}$ are.
\item    
    If $(Y,n,\GA)$ is a common quotient of two systems $(X_i,m_i,\GA)$ ($i=1,2$), 
    then one can form a \emph{relatively independent} joining of $X_1$ with $X_2$ over $Y$
    by taking the measure $\bar{m}$ on $X_1\times X_2$ to be
    \[
    	\bar{m}=\int_Y \mu^{(1)}_y\otimes \mu^{(2)}_y\,dn(y)
    \]
    where $m_i=c_i\cdot\int_Y \mu_y^{(i)}\,dn(y)$, $i=1,2$, are the disintegration with respect
    to the projections.
\item
   An isomorphism (or a $\ta$-twisted isomorphism) $T:(X_1,m_1)\to(X_2,m_2)$ gives rise 
   to the joining ($\ta$-twisted joining) 
   \[
   	\bar{m}=\int_{X_1} \de_x\otimes\de_{T(x)}\,dm_1(x)
   \]
   which is supported on the graph of $T$.
   In particular any (non trivial) element $T$ of the \emph{centralizer} of $(X,m,\GA)$ defines a 
   (non trivial) self joining of $(X,m,\GA)$. 
\item 
   Any ($\tau$-twisted) measurable quotient $p:(X,m,\GA)\to (Y,n,\GA)$ gives rise to the 
   \emph{relatively independent} self joining $(X\times X,\bar{m})$ of $(X,m,\GA)$ given by
   \[
   	\bar{m}=\int_Y \mu_y\otimes\mu_y\,dn(y)
   \]
   where $m=c\cdot \int_Y \mu_y\,dn(y)$ is the disintegration of $m$ with respect to $n$ into 
   a measurable family $Y\to\sP(X)$, $y\mapsto\mu_y$, of probability measures 
   ($\mu_y(p^{-1}\{y\})=1$ and $\mu_{\ga y}=\ga_* \mu_y$ for $n$-a.e. $y\in Y$). 
   This joining need not be ergodic, but can be decomposed into ergodic ones.
\end{enumerate}
Thus understanding joinings between different systems and self joinings of a given system
gives information on isomorphisms between systems, quotients and centralizers.

\medskip

\subsection{An auxiliary Lemma}

We shall need the following technical
\begin{lem}[Pushforward of singular measures] 
\label{L:pushforward}
Let $(X,\mu)$, $(Y,\nu)$ be measure spaces, $Z$ -- a standard Borel space,
$\rh:Y\to Z$ and $x\in X\mapsto \al_x\in\sP(Y)$ be measurable maps
so that
\[
	\int_X \al_x(B)\,d\mu(x)=0\qquad{\rm whenever}\qquad \nu(B)=0
\]
Define a measurable map $X\to\sP(Z)$ by $x\in X\mapsto \be_x=\rh_*\al_x\in\sP(Z)$.
Then 
\begin{enumerate}
\item
The map $\be:X\to\sP(Z)$ is well defined 
in terms of $\al:X\to\sP(Y)$
and $\rh:Y\to\sP(Z)$, all up to null sets.
More precisely if $x\in X\mapsto \al^\prime_x\in\sP(Y)$ agrees $\mu$-a.e. with $\al_x$,
and $\rh^\prime:Y\to Z$ agrees $\nu$-a.e. with $\rh$, then the map
$x\in X\mapsto \be^\prime_x=\rh^\prime_*\al^\prime_x\in\sP(Z)$ agrees with $\be_x$ for
$\mu$-a.e. $x\in X$.
\item
If a countable group $\Gamma$ acts measurably on $X$, $Y$, $Z$ preserving the measure
class of $\mu$ on $X$ and of $\nu$ on $Y$ and such that
\[
	\al_{\ga x}=\ga_* \al_x,\qquad \rh(\ga y)=\ga \rh(y)
\] 
for $\mu$-a.e. $x\in X$, $\nu$-a.e. $y\in Y$ and all $\ga\in\GA$, then
\[
	\be_{\ga x}=\ga_* \be_x
\]
for $\mu$-a.e. $x\in X$ and all $\ga\in\GA$.
\end{enumerate}
\end{lem}
\begin{proof}
(1) For $\mu$-a.e. equality $\be_x=\be^\prime_x$ it suffices to show that for each $E\in\sB(Z)$ 
\[
	\mu\Set{x\in X}{\be_x(E)\neq\be^\prime_x(E)}=0
\]
because $\sB(Z)$ is countably generated ($Z$ is a standard Borel space).
Let $F=\rh^{-1}(E)$, $F^\prime={\rh^\prime}^{-1}(E)\in\sB(Y)$.
We have $\nu(F\vartriangle F^\prime)=0$ and therefore
\[
	\int_X \al_x(F\vartriangle F^\prime)\,d\mu(x)=0
\]
By Fubini $\al_x(F)=\al_x(F^\prime)$ for $\mu$-a.e. $x\in X$. 
At the same time  $\mu$-a.e. $\al_x(F^\prime)=\al^\prime_x(F^\prime)$
and so 
\[
	\be_x(E)=\al_x(F)=\al^\prime_x(F^\prime)=\be^\prime_x(E)
\] 
for $\mu$-a.e. $x\in X$.

(2) For each $\ga$ we have a.e. identities $\al_{\ga x}=\ga_*\al_x$ and $\rh\circ \ga=\ga\circ\rh$
which give rise to 
\[
	\be_{\ga x}=\rh_* (\ga_*\al_x)=(\rh\circ\ga)_*\al_x=(\ga\circ \rh)_*\al_x=\ga_*\be_x
\] 
justified by part (1). 
\end{proof}

\section{Principle bundles with Alignment properties} \label{S:generalCJQ}
The following notion, which we call the \myemph{alignment property}, will play a key 
role in the proofs of our results. 
Section~\ref{S:alignment} contains examples of the alignment property.
Here we shall give the definition, basic properties and a motivating application of this notion.
\begin{defn}\label{D:alignment}
Let $\GA$ be a group, $(X,m)$ be a measure space with a measure class preserving
$\GA$-action, $B$ - a tolpological space with a continuous $\GA$-action,
and $\pi:X\to B$ be a measurable $\GA$-equivariant map.
We shall say that $\pi$ has the \myemph{alignment property} with respect to the
$\GA$-action if $x\mapsto\de_{\pi(x)}$ is the only, modulo $m$-null sets, $\GA$-equivariant 
measurable map from $(X,m)$ to the space $\sP(B)$ of all regular Borel
probability measures on $B$. 
In this case we shall also say that the system 
\[
	(\pi:(X,m)\to B;\GA)
\]
is an \myemph{alignment system}.
\end{defn}
The alignment property depends only on the measure class $[m]$ on a Borel space $X$.
However in all the examples of the alignment phenomena in this paper 
$X$ is a topological space with a continuous $\GA$-action, the map
$\pi:X\to B$ is continuous and $m$ is an infinite $\GA$-invariant measure on $X$.
\bigskip

We start with a list of easy but useful observations about the alignment property.

\medskip

Given a measurable map $\pi:(X,m)\to B$ on a Lebesgue space there is
a well defined measure class $[\nu]$ on $B$ -- 
the ``projection'' $[\nu]=\pi_{*}[m]$  of the measure class $[m]$ on $X$.
It can be defined by taking the projection $\nu=\pi_{*}\mu$ of some \emph{finite} 
measure $\mu$ equivalent to $m$ (being Lebesgue $m$ is $\sigma$-finite).
The measure class $[\nu]$ depends only on $[m]$, and
\[
	\nu(E)=0\quad \Longleftrightarrow\quad m(\pi^{-1}E)=0.
\]
\begin{lem}[Uniqueness] \label{L:align-uniq}
  Let $(\pi:(X,m)\to B;\GA)$ be an alignment system, 
  and $\nu$ a probability measure on $B$ with $[\nu]=\pi_{*}[m]$.
  Then 
  \begin{enumerate}
  \item
  	$\pi$ is the unique, up to $m$-null sets, measurable $\GA$-equivariant map $X\to B$,
  \item
  $b\mapsto \delta_{\pi(b)}\in\sP(B)$ is the unique, up to $\nu$-null sets, 
  measurable $\GA$-equivariant map $B\to\sP(B)$,
  \item
  the identity map is the unique, up to $\nu$-null sets, measurable 
  $\GA$-equivariant map $B\to B$.
  \end{enumerate}
\end{lem}
\begin{proof}
Evident from the definitions.
\end{proof}

\begin{lem}[Conservative]
\label{L:align-cons}
  If $(\pi:(X,m)\to B;\GA)$ is an alignment system and $B$ has more than one point,
  then $\GA$-action on $(X,m)$ has to be conservative. 
\end{lem}
\begin{proof}
  Indeed, otherwise 
  there exists a Borel subset $E\subset X$ with $m(E)>0$ so that 
  $m(\ga_1 E\cap \ga_2 E)=0$ whenever $\ga_1\neq\ga_2\in\GA$. 
  Choose an arbitrary measurable map $p:E\to B$ with $p(x)\neq\pi(x)$,
  extend it in a $\GA$-equivariant way to $\GA E=\bigcup \ga E$ and let $p(x)=\pi(x)$ for
  $x\in X\setminus \GA E$. 
  Then $p$ is measurable $\GA$-equivariant map which does not agree with $\pi$ 
  on a positive measure set $\GA E$ contradicting the alignment property. 
\end{proof}

\begin{lem}[Intermediate Quotients] \label{L:align-interm}
  Let $(X,m,\GA)$ and $(X_0,m_0,\GA)$ be some measure class preserving measurable 
  $\GA$-actions, $p:(X,m)\to (X_0,m_0)$ a $\GA$-equivariant measurable map, 
  $B$ - a topological space with a continuous $\GA$-action, 
  and $\pi_0:(X_0,m_0)\to B$ be a measurable $\GA$-equivariant map.
  If the composition map
  \[
  	\pi:(X,m)\overto{p}(X_0,m_0)\overto{\pi_0} B
  \] 
  has the alignment property, then so does $\pi_0:(X_0,m_0)\to B$. 
\end{lem}
\begin{proof}
Follows from the definitions.
\end{proof}

\begin{lem}[Compact Extensions] \label{L:align-comp-ext}
  Let $p:(X_1.m_1,\GA)\to (X,m,\GA)$ be a compact group
  extension of a $\GA$-action on $(X,m)$, i.e. a compact group $K$ acts on $(X_1,m_1)$
  by measure preserving transformations commuting with the $\GA$-action 
  so that $(X,m)=(X_1,m_1)/K$ with $p$ being the projection.
  If $(\pi:(X,m)\to B;\GA)$ is an alignment system, then the map
  \[
  	\pi_1:(X_1,m_1)\overto{p}(X,m)\overto{\pi} B
  \]
  has the alignment property too.
\end{lem}
\begin{proof}
  If $y\in X_1\mapsto \nu_y\in\sP(B)$ is a measurable $\GA$-equivariant map,
  then 
  \[
  	x\in X\mapsto \mu_x,\qquad \mu_{p(y)}=\int_K \nu_{k y}\,dk
  \] 
  is a $\GA$-equivariant map, and by the alignment property $\mu_x=\de_{\pi(x)}$ for 
  $m$-a.e. $x\in X$. Since Dirac measures are extremal points of $\sP(B)$ it follows
  that $\nu_y=\de_{\pi(p(y))}$ for $m_1$-a.e. $y\in X_1$.
\end{proof}

\begin{lem}[Finite Index Tolerance] \label{L:align-fin-ind}
  Let $(X,m,\Gamma)$ be a measure class preserving action of a countable group $\Gamma$,
  $B$ -- a topological space with a continuous $\Gamma$-action, and $\pi:X\to B$ a measurable
  $\Gamma$-equivariant map. Let $\Gamma^{\prime}<\Gamma$ be a finite index subgroup.
  
  Then $(\pi:(X,m)\to B;\Gamma^{\prime})$ is an alignment system
  iff  $(\pi:(X,m)\to B;\Gamma)$ is also an alignment system.
\end{lem}
\begin{proof}
  Let $\ga_i$, $i=1\dots k$ be the representatives of $\GA^\prime$ cosets.
  Suppose $x\in X\mapsto \nu_x\in\sP(B)$ is $\GA$-equivariant measurable map then 
  $x\mapsto \mu_x=k^{-1} \sum_1^k  \nu_{\ga_i x}$ is $\GA^\prime$-equivariant
  map $X\to\sP(B)$, and therefore is $\mu_x=\de_{\pi(x)}$.
  The fact that $\de_{\pi(x)}$ are extremal points of $\sP(B)$ implies that
  $\nu_x=\de_{\pi(x)}$ for $m$-a.e. $x\in X$.
\end{proof}

\begin{lem}[Products] \label{L:align-products}
  Direct product of alignment systems is an 
  alignment system.
\end{lem}
\begin{proof}
  For $i=1,2$ let $(\pi_i:(X_i,m_i)\to B_i,\GA_{i})$ be alignment systems.
  Set $(X,m)=(X_1\times X_2,m_1\otimes m_2)$, $B=B_1\times B_2$,
  and $\Gamma=\Gamma_{1}\times\Gamma_{2}$. We want to show that
  $\pi(x,y)=(\pi_{1}(x),\pi_{2}(y))$ has the alignment property.

  Suppose that $(x,y)\in X\mapsto \mu_{x,y}\in\sP(B_{1}\times B_{2})$
  is a measurable $\Gamma$-equivariant map. 
  Choose probability measures $m_{i}^{\prime}$ in the measure classes of $m_{i}$.
  Define measurable maps $\nu^{(i)}:(X_{i},m_{i})\to\sP(B_{i})$ by
  \[
  	\nu^{(1)}_{x}(E)=\int \mu_{x,y}(E\times B_{2})\,dm_{2}^{\prime}(y),\qquad
	\nu^{(2)}_{y}(F)=\int \mu_{x,y}(B_{1}\times F)\,dm_{1}^{\prime}(x).
  \] 
  Then $\nu^{(i)}:(X_{i},m_{i})\to \sP(B_{i})$ are measurable and $\Gamma_{i}$-equivariant.
  Thus $\nu^{(i)}_{x}=\delta_{\pi_{i}(x)}$ and since these are extremal points 
  we conclude that a.e.
  \[
  	\mu_{x,y}(E\times B_{2})=\delta_{\pi_{1}(x)}(E),\qquad
	\mu_{x,y}(B_{1}\times F)=\delta_{\pi_{2}(y)}(F).
  \]
  This readily gives $\mu_{x,y}=\delta_{\pi_{1}(x)}\otimes \delta_{\pi_{2}(y)}$.
\end{proof}

\medskip

\subsection{Principle Bundles}
We shall be interested in examples where $\pi:X\to B$ is a principle bundle,
by which we mean that $\pi:X\to B$ is a surjective continuous map between topological spaces,
$L$ is a locally compact group acting continuously and freely on $X$
so that the $L$-orbits are precisely the fibers of $\pi:X\to B$.
In this situation we shall say that $\pi:X\to B$ is a \myemph{principle $L$-bundle.}
An automorphism of a principle $L$-bundle is a homeomorphism 
of $X$ which commutes with the $L$-action, and therefore descends to 
a homeomorphism of $B=X/L$.
\begin{defn}
If $\pi:X\to B$ is a principle $L$-bundle, $\GA$ a group of bundle automorphisms,
$m$ a measure on $X$ so that both $\GA$ and $L$ act on $(X,m)$ by
measure class preserving transformations, and $\pi:(X,m)\to B$
has the alignment property with respect to the $\GA$-action, we shall say
that $(\pi:(X,m)\to B; \GA)$ is a \emph{principle $L$-bundle alignment system},
or an alignment system which is a principle $L$-bundle.
\end{defn}
More specifically we shall consider the following examples:
\begin{exmp}[Homogeneous Spaces]  \label{E:GHGQ}
Let $G$ be a locally compact group, and
$H\normal Q$ - closed subgroups. Set
\[
	X=G/H,\quad B=G/Q, \quad \pi:X\to B,\quad \pi(gH)=gQ.
\]
Observe that $Q$ acts on $G/H$ from the right by $q: gH\mapsto gq^{-1} H$.
This action is transitive on the $\pi$-fibers with $H$ being the stabilizer of every point $gH$.
Thus $L=Q/H$ acts freely on $X$ producing the $\pi$-fibers as its orbits.
Thus $\pi:G/H\to G/Q$ is a principle $L$-bundle.
In this setup the group $G$ and its subgroups act
by bundle automorphisms.

In section \ref{S:alignment} we shall establish the alignment property for such 
bundles under very mild assumptions on $H$ in a semi-simple $G$.
\end{exmp}
\begin{exmp}[Space of Horospheres] \label{E:Hor}
Let $N$ be a complete simply connected Riemannian manifold of pinched negative curvature, 
and let $\pa N$ denote the boundary of $N$. 
For $p,q\in N$ and $\xi\in\pa N$ the Busemann function is defined as
\begin{equation}\label{e:Busemann}
	\be_{\xi}(p,q)=\lim_{z\to\xi} \left[d(p,z)-d(q,z)\right].
\end{equation}
The \myemph{horospheres} in $N$ are the level sets of the Busemann function:
\[
	\hor_{\xi}(t)=\Set{p\in N}{\be_{\xi}(p,o)=t}
\]
where $o\in N$ is some reference point. 
We denote by 
\[
	\Hor(N)=\Set{\hor_{\xi}(t)}{\xi\in\pa N,\ t\in\bR}
\] 
the space of horospheres. 

$\Hor(N)$ fibers over $\pa N$ via $\hor_{\xi}(t)\mapsto \xi$.
This is a principle $\bR$-bundle over $\pa N$, where $\bR$ acts by 
$s:\hor_{\xi}(t)\mapsto \hor_{\xi}(t+s)$.
Since $\ \be_{\xi}(p,q)+\be_{\xi}(q,0)=\be_{\xi}(p,o)\ $ different choices of $o\in N$ change only
the trivialization of this bundle.

The group of isometries of $N$ acts also on the boundary $B=\pa N$
and on the space ${\rm Hor}(N)$ of horospheres because 
$\be_{\ga \xi}(\ga p, \ga o)=\be_{\xi}(p,o)$ for $\ga\in \Isom(N)$.
In the above parametrization this action takes the form
\[
	\ga: \hor_{\xi}(t)\mapsto \hor_{\ga \xi}(t+c(\ga,\xi)),\rmwhere c(\ga,\xi)=\be_{\xi}(\ga o,o).
\]
Note that $c:\Gamma\times \pa N\to \bR$ is an additive cocycle, that is 
\[
	(\gamma^{\prime}\gamma,\xi)=c(\gamma^{\prime},\gamma\cdot \xi)
	+c(\gamma,\xi).
\]
We shall see  in section~\ref{S:alignment} that $\pi: {\rm Hor}(N)\to\pa N$ has the alignment 
property with respect to discrete subgroups $\GA<\Isom (N)$ and $\GA$-invariant
conservative measures (Theorem~\ref{T:Hor}).
\end{exmp}

\medskip

We complete this section by showing how to compute the measurable centralizer,
quotients and self-joinings of any ${\rm II}_{\infty}$ system $(X,m,\Gamma)$ 
admitting a structure of a principle bundle with an alignment property relative
to its base.
\begin{thm}[Centralizers, Self Joinings, Quotients]   \label{T:generalCJQ}
Let $(X,m,\GA)$ be an ergodic infinite measure preserving system 
and $(\pi:(X,m)\to B;\GA)$ a principle  $L$-bundle with the alignment property.
Then the measurable centralizer, self joinings and quotients of $(X,m,\GA)$ can be described as follows:
\begin{enumerate}
\item
	Let ${T}:X\to X$ be some Borel measurable $\GA$-equivariant map.
	Then ${T}(x)=\la_0 x$ for some fixed $\la_0\in L$ and $m$-a.e. $x\in X$.
\item
	Any ergodic self joinings of the $\GA$-action on $(X,m)$ 
	is given by the measure $c\cdot m_{\la_0}$ on $X\times X$ for some $0<c<\infty$ and $\la_0\in L$,
	where
	\[
		m_{\la_0}=\int_{X} \de_x \otimes \de_{\la_0 x}\,dm(x)
	\]
\item
	The only measurably proper $\GA$-equivariant quotients of $(X,m)$
	are of the form $(X,m)/K$ where $K< L$ is a compact subgroup.
\end{enumerate}
\end{thm}
\begin{proof}

\medskip

(1). 
Let ${T}:X\to X$ be an arbitrary Borel measurable $\GA$-equivaraint map.
The map 
\[
	X\overto{{T}} X\overto{\pi} B
\]
is Borel and $\GA$-equivariant as well.
In view of the alignment property of $\pi$ we have
$\pi({T}(x))=\pi(x)$ for $m$-a.e. $x\in X$. 
This allows to define a measurable function $\la:X\to L$ by ${T}(x)=\la_x x$.
We have for $m$-a.e. $x\in X$ and all $\ga\in\GA$:
\[
	\ga \la_{\ga x} x=\la_{\ga x} \ga x={T}(\ga x)=\ga {T}(x)=\ga \la_x x
\]
which, in view of the freeness of the $L$-action, means $\la_{\ga x}=\la_x$ $m$-a.e.
Thus, ergodicitiy of the $\GA$-action on $(X,m)$, implies that $\la_x$ is $m$-a.e.
a constant $\la_0\in L$. 
Hence ${T}(x)=\la_0 x$ for $m$-a.e. $x\in X$ as claimed.

\medskip

(2). 
Let $(X \times X, \bar{m})$ be an ergodic self joining of the $\GA$-action on $(X,m)$.
Then the measure $\bar{m}$ can be disintegrated with respect to its projections on both 
factors as:
\[
	\bar{m}= c_1\cdot \int_{X} \de_{x}\otimes \mu_x\,dm(x)
	   = c_2\cdot \int_{X} \nu_x\otimes \de_x\,dm(x)
\]
with $0<c_1,c_2<\infty$ and $\{\mu_x\}$ and $\{\nu_x\}$ measurable families
of probability measures on $X$, indexed by $x\in X$, satisfying $m$-a.e.
\[
	\mu_{\ga x}=\ga_*\mu_x,\qquad \nu_{\ga x}=\ga_*\nu_x
	\qquad(\gamma\in\Gamma).
\] 
By Lemma~\ref{L:pushforward} we can define a measurable $\Gamma$-equivariant map
\[
	X\ \overto{\mu_{\cdot}}\ \sP(X)\ \overto{\pi_*}\ \sP(B),
	\qquad\text{by}\qquad
	x \mapsto \nu_x \mapsto \pi_*\nu_x.
\]
By the alignment property for $m$-a.e. $x$ the measure $\mu_x$ is supported on 
$\pi^{-1}(\pi(x))=L x$.
Thus $\la x\in L x\mapsto \la\in L$ maps $\mu_x$ to a probability measure $\si_x$ on $L$.
Note that $\{\si_x\}$, $x\in X$, is a measurable family of probability measures on $L$,
which satisfies for $m$-a.e. $x\in X$, every $\ga\in\GA$ and Borel set $E$ on $L$:
\bqrn
	\si_{\ga x}(E)&=&\mu_{\ga x}\left(\Set{\la \ga x}{\la\in E}\right)
		= \ga_*\mu_x \left(\Set{\ga \la x}{\la\in E}\right)\\
	&=& \mu_x\left(\Set{\la x}{\la\in E}\right) =\si_x(E).
\eqrn
Ergodicity of the $\GA$-action on $(X,m)$ implies that 
$\si_x$ is a.e. equal to a fixed probability measure $\si$ on $L$.

The fact that a.e. $\pi_*\mu_x=\de_{\pi(x)}$ means that the measure $\bar{m}$ 
is supported on the set
\begin{equation}\label{e:Relation}
	\sF=\Set{(x,y)\in X\times Y}{\pi(x)=\pi(y)}
\end{equation}
Given a Borel set $E\subs L$ define $F_E=\Set{(x,\la x)}{x\in X,\ \la\in E}\subs \sF$
and observe that
\begin{itemize}
\item 
	$F_E$ is invariant under the diagonal $\GA$-action.
\item
	$\bar{m}(F_E)=0$ if and only if $\si(E)=0$.
\item
	$\sF\setminus F_E=F_{L\setminus E}$.
\end{itemize}
If $\bar{m}$ is ergodic with respect to the diagonal $\GA$-action, for each measurable $E\subset L$
either $\bar{m}(F_E)=\si(E)=0$ or $\bar{m}(F_{L\setminus E})=\si(L\setminus E)=0$.
This implies that $\si$ is a Dirac measure $\de_{\la_0}$ at some $\la_0\in L$, and consequently
$\bar{m}$ is supported on $\Set{(x,\la_0 x)}{x\in X}$ as in the statement of the Theorem.

\medskip

(3).
Let $p:(X,m)\to (Y,n)$ be a $\GA$-equivariant measurably proper quotient. 
Then $m$ can be disintegrated with respect to the quotient map
\[
	m=\int_Y \mu_y\,dn(y)
\] 
where $\mu_y\in\sP(X)$ and $\mu_y(p^{-1}(\{y\}))=1$ for $n$-a.e. $y\in Y$.
Consider the \emph{independent joining relative} to $p:X\to Y$, defined by the measure
\begin{equation}\label{e:rel-ind-join}
	\bar{m}=\int_Y \mu_y\otimes \mu_y\,dn(y)
\end{equation}
on $X\times X$. 
It is $\GA$-invariant and projects onto $m$ in both factors.
Thus the disintegration of $\bar{m}$ into $\GA$-ergodic components
consists of ergodic self joinings of $(X,m,\GA)$ i.e.
\[
	\bar{m}=\int_L m_\la\,d\si(\la),\rmwhere m_\la=\int_X \de_x\otimes\de_{\la x}\,dm(x)
\]
and $\si$ is a probability measure on $L$.
In particular $m$ is supported on the set 
\[
	\sF=\Set{(x,x^\prime)}{\pi(x)=\pi(x^\prime)}=\Set{(x,\la x)}{x\in X,\ \la\in L}\subset X\times X.
\]
A comparison with (\ref{e:rel-ind-join}) yields that for $n$-a.e. $y\in Y$ the 
measure $\mu_y$ is supported on a single $L$-orbit and moreover
for $\mu_y$-a.e. $x$ 
\[
	\int_L f(\la x)\,d\si(\la) = \int_X f(x^\prime)\, d\mu_y(x^\prime),\qquad (f\in C_c(X)).
\]
Since the roles of $x$ and $x^\prime$ are symmetric, $\si$ is a symmetric measure
i.e. $d\si(\la)=d\si(\la^{-1})$. 
Moreover, for any $f\in C_c(X)$ and $\mu_y$-a.e. $x$ we have 
\bqrn
	&& \int_L\int_L f(\la_2\la_1 x)\,d\si(\la_1)\,d\si(\la_2)
	     =\int_L \left(\int_X f(\la_2 x_1)\,d\mu_y(x_1)\right)\,d\si(\la_2)\\
		&&=\int_X\left(\int_L f(\la_2 x_1)\,d\si(\la_2)\right)\,d\mu_y(x_1)		=\int_X f\,d\mu_y=\int_L f(\la x)\,d\si(\la)
\eqrn
which implies that $\si*\si=\si$.
\begin{lem}\label{L:idempotent-compact}
A symmetric probability measure $\si$ on a locally compact group $L$ 
satisfies $\si*\si=\si$ if and only if $\si=m_K$ -- Haar measure on a compact subgroup $K< L$.
\end{lem}
\begin{proof}
The "if" part is evident. 
For the "only if" part assume $\si$ is symmetric and $\si*\si=\si$ and
let $K=\supp(\si)$. 
$K$ is a closed subgroup of $L$. 
Indeed $K^{-1}=K$ and $K\cdot K\subs K$. 
To see the latter, given $k_1,k_2\in\supp(\si)$ and a neighborhood $U$ of $k_1\cdot k_2$, 
coose open neighborhoods $V_i$ of $k_i$ so that $V_1\cdot V_2\subset U$,
and note that 
\[
      \si(U)=\si*\si(U)\ge \si(V_1)\cdot\si(V_2)>0.
\]
As $U$ was arbitrary, it follows that $k_1\cdot k_2\in K$.

Now let $P_\si$  be the Markov operator
\[
	(P_\si f)(k)=\int_K f(k k^\prime)\,d\si(k^\prime).
\]
which is defined on $C_0(K,\bR)$ and takes values in $C_0(K,\bR)$.
It is a projection because $P_\si^2=P_{\si*\si}=P_\si$.
If $g\in C_0(K,\bR)$ is a $P_\si$-invariant function, then
the closed set $A_g=\Set{k\in K}{g(k)=\max g}$ satisfies
$A_gk^\prime=A_g$ for $\si$-a.e. $k^\prime\in K$, and $K=\supp(\si)$ yields 
$A_g=K$ and so $g=const$. 
This implies that $K$ is compact for $C_0(K,\bR)$ contains a non-trivial constant function.
Hence for $f\in C(K)$ and $k_0\in K$ we have
\[
	\int_K f(k)\,d\si(k) = (P_\si f) (e)=(P_\si f)(k_0)=\int_K f(k_0k)\,d\si(k)
\]
which means that the probability measure $\si$ is left invariant on the compact group $K$.
\end{proof}

\begin{rem}
In the above Lemma the assumption that $\si$ is symmetric is redundant:  
we used only $\si=\si*\si$ to show that $K=\supp(\si)$ is a closed sub semigroup
of $L$ which is compact; but  compact sub semigroups in topological
groups are known to form subgroups.
\end{rem}

Returning to the description of the quotient $(Y,n)$ of $(X,m)$,
we observe that for $n$-a.e. $y\in Y$ the measure $\nu_y$ is 
supported and equidistributed on a single $K$-orbit in $X$,
so that $Y$ can be identified with $X/K$ and $n$ with $m/K$.
Note also that the $\GA$-action on $X$ descends to an action on $X/K$
because $\GA$ and $K<L$ commute. 
This completes the proof of Theorem~\ref{T:generalCJQ}.
\end{proof}
 
\section{Main Examples of the Alignment Property} \label{S:alignment}
We start with the geometric setup (recall Example~\ref{E:Hor}):
\begin{thm}[Space of Horospheres]
\label{T:Hor}
Let $N$ be a complete simply connected negatively curved Riemannian manifold,
$X={\rm Hor}(N)$ -- the space of horoshperes, $B=\pa N$ -- the boundary, $\pi:X\to B$
the natural projection. Let $\GA<\Isom(N)$ be a discrete group, and $m$ be some
Borel regular $\GA$-invariant measure with full support on $X$.

Then $(\pi:(X,m)\to B;\Gamma)$ is the alignment system if and only if 
the $\GA$-action on $(X,m)$ is conservative.
\end{thm}

\begin{proof}
Lemma~\ref{L:align-cons} provides the "only if" direction.
The content of this Theorem is the "if" direction. 
We assume that the $\GA$-action on $(X,m)$ is conservative
and $x\in X\mapsto\mu_x\in\sP(B)$ is a measurable $\GA$-equivariant map,
which should be proven to coincide $m$-a.e. with $\de_{\pi(x)}$.
Assuming the contrary, the set
\[
	A=\Set{x\in X}{\mu_x(\{\pi(x)\})<1}\qquad{\rm has}\qquad m(A)>0.
\]
Note that $A$ is $\GA$-invariant and we may assume that $\mu_{\ga x}=\ga_*\mu_x$
for \emph{all} $x\in A$ and all $\ga\in\GA$. 
For $x\in A$ denote by $\nu_x$ the normalized restriction of $\mu_x$ to $B\setminus\{\pi(x)\}$. 
Then also $\{\nu_x\}$, $x\in A$, is $\GA$-equivariant: $\nu_{\ga x}=\ga_*\nu_x$.
 
Let $\rh$ be some metric on $B=\pa N$, e.g. the visual metric from the base $o\in N$,
and let $U_{\xi,\eps}=\Set{\et\in B}{\rh(\xi,\et)<\eps}$, 
$K_{x,\eps}=B\setminus U_{\pi(x),\eps}$.
Set
\[
	A_\eps=\Set{x\in A}{\nu_x(K_{x,\eps})>1/2}
\]
As $A=\bigcup_1^\infty A_{1/k}$, there exists $\eps>0$ so that $m(A_\eps)>0$.
By Luzin's theorem there exists a (compact) subset $C\subs A_\eps$
with $m(C)>0$ so that the map $x\in C\mapsto \nu_x\in\sP(B)$ is \emph{continuous} on $C$.
Since $\GA$ is conservative on $(X,m)$ and $m$ is positive on non-empty open sets, 
for $m$-a.e. $x\in C$ there exists an infinite sequence of elements $\ga_n\in \GA$ so that
\[
	\ga_n x\to x, \rmand \ga_n x\in C\subs A_\eps,\ n\in\bN.
\] 
Let us fix such an $x_0$ and the corresponding infinite sequence $\{\ga_n\}$.
Let $X=K_{x_0,\eps}$ - a compact subset of $B=\pa N$ with $\rh(\pi(x_0),K)\ge\eps>0$.
We shall show that for large $n$ the set $\ga_n^{-1} K$ and $K$ are disjoint.
Then 
\[
	\nu_x(K)+\nu_x(\ga_n^{-1}K)=\nu_x(K\cup \ga_n^{-1} K)\le 1
\]
for large $n$, and we shall arrive at a contradiction because
\[
	1/2<\nu_x(K)=\lm \nu_{\ga_n x}(K)=\lm (\ga_n)_*\nu_x(K)=\lm \nu_x(\ga_n^{-1} K).
\] 
In order to show that $\ga_n^{-1}K$ and $K$ are disjoint for large $n$,
we look at the unit tangent bundle $SN$ of $N$. 
It is homeomorphic to
\[
	\Set{(\et,\xi,t)}{\et\neq\xi}\quad{\rm via}\quad v\in SN\mapsto (v^+,v^-,t(v))
\]
were $v^+\neq v^-\in\pa N$ denote the forward and backward end points
of the geodesic in $N$ determined by $v$, and $t(v)=\be_{v^{+}}(p(v),o)\in\bR$,
where $p(v)\in N$ is the base point of the unit tangent vector $v$.

Let $\xi_0=\pi(x_0)$.
The horoshere $x_0=\hor_{\xi_0}(t_0)$ corresponds in one-to-one fashion to 
$\Set{(\xi_0,\et,t_0)}{\et\neq\xi_0}$
-- the set of unit vectors based at points of $\hor_{\xi_0}(t_0)$
and pointing towards $\xi_0\in\pa N$ (this is the usual identification between the horosphere
as a subset of $N$ and as the stable leaf in $SN$).
The assumption that $\ga_n x_0\to x_0$ in $X$ means that 
$\ga\xi_0\to \xi_0$ and $t_n\to t_0$ in $\bR$,
where $\ga_n :\hor_{\xi_0}(t_0)\mapsto \hor_{\ga_n \xi_0}(t_n)$.

Let $n_0\in\bN$ be such that $\rh(\ga_n\xi_0,\xi_0)<\eps/2$ and $|t_n-t_0|<1$
for all $n\ge n_0$.
The set 
\[
	Q=\Set{u\in SN}{\rh(u^+,u^-)\ge \eps/2,\ t(u)\in[t_0-1,t_0+1]}
\]
is compact in $SN$. 
$\Isom(N)$ acts properly discontinuously on $N$ and on $SN$. 
Thus the infinite sequence $\{\ga_n\}$ in the discrete subgroup $\GA<\Isom(N)$
eventually moves $Q$ away from itself. 
In other words, there exists $n_1$ so that
\[
	\ga_n Q\cap Q=\emptyset,\qquad \forall n\ge n_1.
\]
We claim that $\ga_n^{-1}K\cap K=\emptyset$ for all $n\ge \max\{n_0,n_1\}$. 
Indeed suppose that $\et\in \ga_n^{-1} K\cap K$ for some $n\ge n_0$, 
and let $v\in SN$ be (the unique) vector with $(v^+,v^-,t(v))=(\xi_0,\et,t_0)$.
Then $v\in Q$ and $\ga_n v$, corresponding to $(\ga_n \xi_0,\ga_n\et,t_n)$, also lies in $Q$.
This is possible only if $n<n_1$.
This completes the proof of the Theorem.
\end{proof}

\bigskip

For the proof of the alignment property for maps between homogeneous
spaces we shall need the following general result.
\begin{thm}[Homogeneous spaces]\label{T:gen-alignment}
Let $G$ be a locally compact group, 
$H< Q< G$ -- closed subgroups and $\GA< G$ -- a discrete subgroup.
Suppose that there exists an open cover
\[
	G/Q\setminus\{eQ\}=\bigcup V_i
\] 
and closed subgroups $T_i<H$ so that for each $i$:
\begin{enumerate}
\item
	$V_i$ is a $T_i$-invariant set and $T_i$ acts properly on $V_i$,
\item
	the $\GA$-action on $(G/T_i,m_{G/T_i})$ is conservative.
\end{enumerate}
Let $(X,m)=(G/H,m_{G/H})$, $B=G/Q$ and $\pi:X\to B$ be the natural projection.
Then $\pi:(X,m)\to B$ has the alignment property with respect to the $\GA$-action.
\end{thm}
\begin{proof}
Let $X\to\sP(G/Q)$, $x\mapsto \mu_x$, be a fixed Borel map satisfying
\[
	\mu_{\ga x}=\ga_*\mu_x
\] 
for all $\ga\in\GA$ and $m$-a.e. $x\in X$.
We shall prove that $\mu_x=\de_{\pi(x)}$ for $m$-a.e. $x\in X$ by reaching a 
contradiction starting from the assumption that $m(A)>0$ where
\[
	A=\Set{x\in X}{\mu_x(B\setminus\{\pi(x)\})>0}.
\]
Let $x\mapsto g_x$ be some Borel cross section of the projection $G\to G/H$.
As $B\setminus\{\pi(x)\}=g_x(G/Q\setminus \{eQ\})=g_x \left(\bigcup V_i\right)$ 
we have
\[
	A=\Set{x\in X}{\mu_x(\cup_i g_xV_i)>0}\subs
	 \bigcup_i \Set{x\in X}{\mu_x( g_xV_i)>0}.
\]
The set of indices $i$ in the assumption of the theorem can always be taken finite or countable
because all the homogeneous spaces in question are separable (in the applications
below the set is actually finite). 
Thus $m(A)>0$ implies that for some $i$ the set 
\[
	A_i=\Set{x\in X}{\mu_x(g_x V_i)>0}
\]
has positive $m$-measure.
Let $x\in A\mapsto \nu_{x}\in\sP(G/Q)$ denote the normalized restrictions of $\mu_x$ to $g_x V_i$:
\[
	\nu_x(E)=\mu_x(E\cap g_x V_i)/\mu_x(g_x V_i)\qquad 
	(E\subset G/Q \ \text{measurable}).
\]
We still have $\nu_{\ga x}=\ga_* \nu_x$ for a.e. $x\in A_i$ and all $\ga\in\GA$.
Given a compact subset $K$ of $V_i$ let 
\[
	A_{i,K}=\Set{x\in A_i}{\nu_x(g_x K)>1/2}.
\] 
and choose $K\subset V_i$ large enough so that $m(A_{i,K})>0$.

By Luzin's theorem, there exists a (compact) subset 
$C\subs A_{i,K}$ with $m(C)>0$, so that both $g_x\in G$ and $\nu_x\in\sP(G/L)$ 
vary \emph{continuously} on $x\in C$.
We can also assume that $\nu_{\ga x}=\ga_* \nu_x$ for \emph{all} $x\in C$ and all $\ga\in\GA$.

We shall now use the assumption that the $\GA$-action on $(X^\prime,m^\prime)=(G/T_i,m_{G/T_i})$ is conservative
in order to obtain the following
\begin{lem}\label{L:cons}
For $m$-a.e. $x\in C$ there exist sequences: 
$\ga_n\to\infty$ in $\GA$, $u_n,v_n\to e$ in $G$, and $t_n\to\infty$ in $T_i$,
so that:
\[
	\ga_n=g_x v_n t_n u_n^{-1} g_x ^{-1},\qquad \ga_n x\in E, \qquad\ga_n x\to x.
\]
\end{lem}
\begin{proof}
Let $U$ be a neigborhood of $e\in G$. 
Choose a smaller neighborhood $W$ for which $e\in W=W^{-1}$ and $W^2\subset U$.
Since $x\in C\mapsto g_x\in G$ is continuous, the compact set $C\subset X$
can be covered by (finitely many) subsets $C=\bigcup B_j$ of small enough size to ensure
that $g_x^{-1}g_y\in W$ whenever $x$ and $y$ lie in the same $B_j$.
Consider the subsets $E_j\subset X^\prime$ defined by
\[
	E_j=\Set{g_x w T_i}{x\in B_j,\ w\in W\cap H}.
\]
Then $m^\prime(E_j)>0$ whenever $m(B_j)>0$. For each such $j$ apply the following argument:
$\GA$ is conservative on $(X^\prime,m^\prime)$, i.e. $m^\prime$-a.e. point of $E_j$ is recurrent.
This implies that for $m$-a.e. $x\in B_j$ and $m_H$-a.e. $w\in W\cap H$ there exists a non trivial 
$\ga\in\GA$ such that 
$\ga g_xw T_i=g_y w^\prime T_i$ where $y=\ga x\in B_j$ and $w^\prime\in W\cap H$.
Hence for some $t\in T_i$ we have $\ga g_x w=g_y w^\prime t$ and
\[
	\ga=g_x (g_x^{-1}g_y)w^\prime t w^{-1} g_x^{-1}=g_x u t v g_x^{-1}
\]
where $u=(g_x^{-1}g_y)w^\prime\in W^2\subset U$ and $v=w^{-1}\in W^{-1}\subset U$.

This shows (by applying these arguments to all $B_j$ of positive measure) that for $m$-a.e. $x\in C\subset X$
there exists $\ga\in\GA$ with $\ga x\in C$ and $\ga= g_x u t v g_x^{-1}$ where $u,v\in U$ and $t\in T_i$.
By passing to a sequence $\{U_n\}$ of neighborhoods shrinking to identity in $G$, we obtain 
the sequence $\ga_n=g_x v_n t_n u_n^{-1} g_x ^{-1}$ with $u_n,v_n\to e$. 
$\GA$ is discrete in $G$. Thus $\ga_n\to\infty$ in $G$ which yields $t_n\to\infty$ in $T_i$.
By construction $\ga_n x\in E$, and $\ga_n x\to x$ from the above form of $\ga_n$.
\end{proof}

We return to the proof of the Theorem.
Let $K^\prime\subset V_i$ be a compact set containing $K$ in its interior.
For a.e. $x\in C$ let $\ga_n$, $u_n$, $v_n$ and $t_n$ be as in Lemma~\ref{L:cons}.
Choose  a symmetric neighborhood $U$ of $e\in G$ small enough to ensure 
$U K \subset K^{\prime}$.   
For $n$ large enough, $u_n, v_n \in U$. Thus
\bqrn
	\nu_{\ga_n x}(g_x K)&=&\ga_n \nu_x(g_{x} K)= \nu_x(\ga_n^{-1} g_{x} K) \\
	&=& \nu_x(g_x u_n t_n^{-1} v_n^{-1}g_x^{-1}g_{x} K)\le \nu_x(g_x u_n t_n^{-1} K^\prime)
\eqrn
Since $\ga_n x\to x$ in $C$, the left hand side converges to $\nu_x(g_x K)>1/2$
which gives
\[
	\liminf_{n\to\infty} \nu_x(g_x u_n t_n^{-1} K^\prime)>1/2.
\]
At the same time $t_n\to\infty$ in $T_i$ and in view of the properness of the $T_i$-action on $V_i$,
for large $n$, $t_n^{-1} K^\prime$ is disjoint from $UK\subset K^\prime$.
Hence
\[
	g_x u_n t_n^{-1}K^\prime\cap g_x K\subs g_x u_n (t_n^{-1} K^\prime\cap UK)=\emptyset
\]
which yields
\[
	\nu_x(g_x u_nt_n^{-1} K^\prime)\le 1-\nu_x(g_x K)<1/2.
\]
The obtained contradiction completes the proof of Theorem~\ref{T:gen-alignment}.
\end{proof} 

\bigskip

Let us illustrate this general result in two concrete examples, and then in a more general situation.

\begin{cor}  \label{C:F2}
Let $\locfld$ be a local field, $G=\SL_2(\locfld)$, 
and $\GA< G$ a discrete subgroup acting conservatively on $\locfld^2$ w.r.to the Haar measure.
Then the projection $\pi:\locfld^2\setminus\{0\}\to \PrjSp{\locfld}{1}$ has the alignment property with 
respect to the $\GA$-action.
\end{cor}
\begin{proof}
Denote by $H$ the stabilizer in $G$ of $e_1=(1,0)\in\locfld^2$, and by $Q$ the stabilizer in $G$ 
of the projective point $[e_1]=\locfld e_1\in\PrjSp{\locfld}{1}$. 
Then $G/Q\iso\PrjSp{\locfld}{1}$ and $(G/H,m_{G/H})\iso (\locfld^2\setminus\{0\},{\rm Haar})$, 
 the $H$-action on $G/Q$ acts properly discontinuously
(and transitively) on the complement $V=G/Q\setminus \{eQ\}$ of the fixed point
$\{eQ\}$. The assumptions of Theorem~\ref{T:gen-alignment} are satisfied.
\end{proof}

\begin{cor}\label{C:Fn}
Let $\locfld$ be a local field, $G=\SL_n(\locfld)$ and $\GA< G$ be a lattice.
Then the projection $\pi:\locfld^n\setminus\{0\}\to \PrjSp{\locfld}{n-1}$ has the alignment property 
with respect to the $\GA$-action and the Haar measure on $\locfld^n$.
\end{cor}
\begin{proof}
Let $H$ be the stabilizer in $G$ of $e_1=(1,0,\dots,0)\in\locfld^n$, and let $Q$ be the stabilizer of the 
projective point $[e_1]=\locfld e_1\in\PrjSp{\locfld}{n-1}$. 
Then $(G/H,m_{G/H})\iso (\locfld^n\setminus\{0\},{\rm Haar})$ and
$G/Q\iso\PrjSp{\locfld}{1}$. 
For $i=2,\dots,n$ let 
\[
	V_i=\Set{[(x_1,\dots,x_n)]\in\PrjSp{\locfld}{n-1}}{x_i\neq 0},\qquad
	T_i=\Set{I+t E_{1,i}}{t\in\locfld}
\]
where $I$ denotes the identity matrix and $E_{j,k}$ the elementary matrix with $1$ in the $j,k$-place
and zeros elsewhere. This system satisfies the assumptions of Theorem~\ref{T:gen-alignment}.
Indeed the only non-elementary condition is conservativity of the $\GA$-action on $G/T_i$.
By Moore's ergodicity theorem $T_i$ acts ergodically on $G/\GA$, which is equivalent 
to the ergodicity of the $\GA$-action on $G/T_i$.
\end{proof}

The above are particular cases of the following more general: 
\begin{thm}\label{T:GH}
Let $G=\prod {\mathbf G}_{\alpha}(\locfld_{\alpha})$ be a semi-simple group and
$H<G$ be a ``suitable'' subgroup as in Definition~\ref{D:suitableH}, i.e. $H$ is pinched 
$\check{H}\normal H <\hat{H}$ between certain unimodular subgroups 
$\check{H},\hat{H}\normal Q$ associated to a parabolic $Q<G$.
Let $\Gamma<G$ be a lattice (not necessarily irreducible).
Denote $(X,m)=(G/H,m_{G/H})$, $(\check{X},\check{m})=(G/\check{H},m_{G/\check{H}})$, 
$(\hat{X},\hat{m})=(G/\hat{H},m_{G/\hat{H}})$ and $B=G/Q$.
Then
\begin{enumerate}
\item
	the natural projections 
	\[
		\check{\pi}:(\check{X},\check{m})\to B,\quad
		\pi:(X,m)\to B,\quad
		\hat{\pi}:(\hat{X},\hat{m})\to B
	\]
	have the alignment property with respect to the action of $\Gamma$; 
\item
	the systems $\check{\pi}:\check{X}\to B$, $\hat{\pi}:\hat{X}\to B$ 
	are principle bundles with structure groups 
	$\check{L}=Q/\check{H}$ and $\hat{L}=Q/\hat{H}$;
\item
	$\hat{L}$ is Abelian and has no non-trivial compact subgroups, 
	and $\check{L}$ is a compact extension of $\hat{L}$.
\end{enumerate}
\end{thm}
%
\begin{proof}[Proof of Theorem~\ref{T:GH}]
Assertions (2) and (3) just record some properties of the construction \ref{D:suitableH}
which will be used later (cf. Theorem~\ref{T:char-maps}).
So it only remains to prove  the alignment property of the three systems (1).
The natural projections are nested
\[
	\check{X}\to X\to \hat{X}\overto{} B,\qquad 
	g\check{H}\mapsto gH\mapsto g\hat{H}\mapsto gQ.
\]
So by Lemma~\ref{L:align-interm} it suffices to prove the alignment property for  
the $\Gamma$-action on $\check{\pi}:(\check{X},\check{m})\to B$.

The group $G$, as well as $\check{H}\normal Q$, are product groups:
\[
	G=\prod_{\alpha\in A} G_{\alpha},\qquad
	\check{H}=\prod_{\alpha\in A} \check{H}_{\alpha},\qquad
	Q=\prod_{\alpha\in A} Q_{\alpha}
\]
formed by the $\locfld_{\alpha}$-points of the corresponding $\locfld_{\alpha}$-groups.
So $\check{\pi}:\check{X}\to B$ splits as a product of projections 
\[
	\check{\pi}_{\alpha}:\check{X}_{\alpha}=G_{\alpha}/H_{\alpha}\ \overto{} \ 
	B_{\alpha}=G_{\alpha}/Q_{\alpha}.
\]
If $\Gamma$ is reducible, then some subgroup of finite index $\Gamma^{\prime}$
splits as a product of irreducible lattices $\Gamma_{i}<\prod_{\alpha\in A_{i}}G_{\alpha}$,
where $A=\cup_{i} A_{i}$ is some non-trivial partition. 
Thus Lemmas~\ref{L:align-fin-ind} and \ref{L:align-products} allow to transfer
the alignment property from $\Gamma^{(i)}$-actions on 
$\check{X}^{(i)}=\prod_{\alpha\in A_{i}}\check{X}_{\alpha}\to B^{(i)}=\prod_{\alpha\in A_{i}}B_{\alpha}$ to that of 
$\Gamma$ on $\check{X}\to B$.
Therefore we may assume that $\Gamma$ itself is irreducible in $G$.

We now point out two important properties of $\check{H}$: 
\begin{itemize}
\item[(i)]
	$\check{H}$ is generated by unipotent elements in its factors $\check{H}_{\alpha}$.
\item[(ii)]
	If $g\in G$ satisfies $g^{-1}\check{H} g<Q$ 
	then necessarily $g\in Q$.
\end{itemize}
The first property follows from the construction, in fact, $\check{H}_{\alpha}$ is
generated by all the roots of $G_{\alpha}$ contained in $Q_{\alpha}$.
This also explains the second statement. 

Next, let $\{T_i\}$ be a (finite) family of one parameter unipotent 
subgroups (collected from different factors) in $\check{H}$ generating $\check{H}$.
The $T_i$-action on the projective variety $B=G/Q$ is algebraic. 
Hence $G/Q$ decomposes as $F_i\sqcup V_i$ where $F_i$ is the set of $T_i$-fixed
points and $V_i$ is the union of free $T_i$-orbits 
($T_i\iso\locfld$ has no proper algebraic subgroups).
The intersection $\bigcap F_i$ consists of $\check{H}$-fixed points, 
i.e. points $gQ$ such that $\check{H}gQ=gQ$. 
Property (ii) above yields
\[
	\bigcap F_i=\{eQ\},\qquad{\rm so\ that}\qquad
	\bigcup V_i=G/Q\setminus\{eQ\}.
\]
By Moore's ergodicity theorem $T_i$ acts ergodically on $G/\GA$ and so $\GA$-action on $G/T_i$ is ergodic and hence conservative with respect to the Haar measure $m_{G/T_i}$.
The $T_i$-action on $V_i$ is properly discontinuous.
It therefore follows from Theorem~\ref{T:gen-alignment} that $\pi:G/\check{H}\to G/Q$ 
has the alignment property. 
\end{proof}

\section{Joinings of Principle bundle actions}\label{S:boundary}
In this section we consider two ergodic infinite measure actions $(X_i,m_i,\GA)$, $i=1,2$, 
which are either isomorphic or, more generally, admit a joining (which we shall assume to be ergodic). 
Assuming that these actions have alignment systems $\pi_i:(X_i,m_i)\to B_i$
which are principle bundles with \emph{amenable} structure groups we 
relate $\GA$-actions on $B_i$ and the structure groups.
\begin{thm}[Joinings and Boundary maps] \label{T:char-maps}
Let $(X_i,m_i,\GA)$, $i=1,2$, be two ergodic infinite measure
preserving actions of some countable group $\Gamma$, 
which admit an ergodic $\ta$-twisted joining $\bar{m}$ 
where $\ta\in\Aut\GA$.  
Assume that  
\begin{itemize}
\item
	$(\pi_i:(X_i,m_i)\to B_i;\GA)$ are principle $L_i$-bundles with the alignment property;
\item
	$L_i$ are amenable groups and $B_i$ are compact spaces with continuous $\GA$-actions. 
\end{itemize}
Denote by $[\nu_i]={\pi_i}_*[m_i]$ the projection measure class on $B_i$.
Then
\begin{enumerate}
\item 
	There exists a measurable bijection $\ph:(B_1,\nu_1)\to (B_2,\nu_2)$ so that 
	\[ 
		[\ph_*\nu_1]=[\nu_2],\qquad \ph(\ga b)=\ga^\ta \ph(b)
	\] 
	for $\ga\in\GA$ and $\nu_1$-a.e. $b\in B_1$, 
	the map $\ph$ being uniquely defined, up to null sets, by its properties.
	The joining $\bar{m}$ is supported on the set
	\[ 
		\sF_\ph=\Set{(x_1,x_2)\in X_1\times X_2}{\ph(\pi_1(x_1))=\pi_2(x_2)}.
	\] 
\item
	Structure groups $L_i$ contain compact subgroups $K_i$, so that  
	the joining $\bar{m}$ is $K_{1}\times K_{2}$-invariant; 
	the quotient systems $(X_i^\prime,m_i^\prime)=(X_i,m_i)/K_i$
	(see Lemma~\ref{L:quotient-by-comp}) admit a measurable isomorphism
	\[
		T^\prime:(X_1^\prime,m_1^\prime)\to(X_2^\prime,c\cdot m_2^\prime),\qquad
		T^\prime(\ga x_1)=\ga^\ta T^\prime(x_1).
	\]
	The natural projections $X_{i}\to X_{i}^{\prime}$ ($i=1,2$) define 
	a $\Gamma$-equivariant quotient map from $\sF_{\phi}\subset X_{1}\times X_{2}$
	with the joining measure $\bar{m}$ to the graph of $T^{\prime}$ with the measure
	\[
		\bar{m}^{\prime}=const\cdot \int_{X^{\prime}_{1}}
			\delta_{x}\otimes\delta_{T^{\prime}(x)}\,dm^{\prime}_{1}(x).
	\]
\end{enumerate}
Furthermore, if $K_{i}$ are normal in $L_{i}$ (for example if the joining
$\bar{m}$ is a graph of an isomorphism, then $K_{i}=\{e_{i}\}$ are trivial), then
\begin{itemize}
\item[{\rm (3)}]
	The alignment systems $\pi^{\prime}_{i}:X_{i}^{\prime}\overto{}B_{i}$
	are principle $\Lambda_{i}$-bundles where $\Lambda_{i}=L_{i}/K_{i}$;  
	there exists a continuouse isomorphism 
	$\rho:\Lambda_{1}\overto{\iso}\Lambda_{2}$ and the measurable 
	$\tau$-twisted isomorphism $T^{\prime}:X^{\prime}_{1}\to X_{2}^{\prime}$ satisfies:
	\[ 
		T^{\prime} (\la_1 x)=\rh(\la_1) T^\prime(x)
	\] 
	for $m^{\prime}_1$-a.e. $x\in X^\prime_1$ and $m_{\LA_1}$-a.e. $\la_1\in \LA_1$.
\end{itemize} 
\end{thm}

The following commutative diagram schematically summarizes these statements:
\[
     	\begin{CD}
     	  (X_1,m_1)  @>{\bar{m}}>>  (X_2,m_2) \\
      	 @VV{/K_1}V   @VV{/K_2}V \\
   	  (X^\prime_1,m^\prime_1)  @>{T^\prime}>>  (X^\prime_2,m^\prime_2) \\
     	 @V{\pi^\prime_1}V{\circlearrowleft\LA_1}V   @V{\iso\LA_2\circlearrowright}V{\pi^\prime_2}V \\
      	  (B_1,[\nu_1]) @>{{\ph}}>>  (B_2,[\nu_2])
    	 \end{CD}
\]

\begin{proof}[Proof of Theorem~\ref{T:char-maps}]
The joining $\bar{m}$ can be presented as 
\begin{equation}\label{e:barm}
	\bar{m}=c_1\cdot\int_{X_1} \de_x\otimes\mu_x\,dm_1(x)=c_2\cdot\int_{X_2}\nu_y\otimes\de_y\,dm_2(y)
\end{equation}
where the measurable maps $x\in X_1\mapsto\mu_x \in\sP(X_2)$, $y\in X_2\mapsto \nu_y\in\sP(X_1)$ 
satisfy almost everywhere
\[
	\mu_{\ga x}=\ga^\ta_*\mu_x,\qquad \nu_{\ga^\ta y}=\ga_*\nu_y.
\]
The families $\{\mu_x\}$ and $\{\nu_y\}$ are uniquely determined, 
up to null sets, by the joining $\bar{m}$. 
Thus Lemma~\ref{L:pushforward} allows us to define a measurable map $p:X_1\to\sP(B_2)$ 
\[
	p(x)={\pi_2}_*\mu_x,\qquad{\rm satisfying}\qquad p(\ga x)=\ga^\ta_* p(x)
\]
for all $\ga\in\GA$ and $m_1$-a.e. $x\in X_1$.
\begin{clm}
There exists a measurable map $\ps:B_1\to\sP(B_2)$ so that $\ps(\ga b)=\ga^\ta_*\ps(b)$
for all $\ga\in\GA$ and $\nu_1$-a.e. $b\in B_1$.
\end{clm}
The idea is to construct  $\ps$ by "averaging" $p:X_1\to \sP(B_2)$ over the action of the 
amenable group $L_1$ on $X_1$.
\[
     \begin{CD}
      X_1  @>{\mu_\cdot}>>  \sP(X_2) \\
       @A{L_1\circlearrowright}AA   @VV{{\pi_2}_*}V \\
        B_1 @>{{\ps}}>>  \sP(B_2)
     \end{CD}
\]
Indeed, let $b\in B_1\mapsto x_b\in X_1$ be some Borel cross section of $\pi:X_1\to B_1$.
If $\si$ is a probability measure on $L_1$ which is absolutely continuous with respect to $m_{L_1}$,
then for $\nu_1$-a.e. $b\in B_1$ the expression 
\[
	\int_{L_1} p(\la x_b)\,d\si(\la)
\]
is a well defined probability measure on $B_2$.
We shall apply this construction to an asymptotically invariant  sequence $\si_n$ of absolutely continuous 
probability measures on $L_1$, i.e. measures $\si_n\preceq m_{L_1}$ such that 
\[
	\|\la_1\si_n-\si_n\|\to 0
\] 
for any $\la_1\in L_1$. 
Existence of such a sequence $\{\si_n\}$ is one of equivalent definitions of amenability of $L_1$. 
Define a sequence of maps $\ps_n:B_1\to\sP(B_2)$ by
\[
	\ps_n(b) = \int_{L_1} p(\la x_b)\,d\si_n(\la).
\]
We shall now apply a compactness argument and pass to a convergent subsequence $\ps_n\to\ps$.
Indeed the space $\DE$ of (equivalence classes of) measurable maps $B_1\to\sP(B_2)$ 
(where two maps which agree $\nu_1$-a.e. are identified), 
can be viewed as a closed convex subset of the unit ball in the Banach space
\[
	W=L^\infty\left((B_1,\nu_1)\to C(B_2)^*\right)
\] 
which is the dual of the Banach space 
$V=L^1((B_1,\nu_1)\to C(B_2))$. 
Thus, endowed with the weak-* topology, $\DE$ is compact.
Passing to a subsequence, we may assume that $\ps_n\to\ps\in\DE$,
where the convergence is in the weak-* sense in $W=V^*$. 
Passing to a further subsequence we can guarantee that for $\nu_1$-a.e. $b\in B_1$
$\ps_n(b)\to \ps(b)$ in weak-* sense in $\sP(B_2)$.

The chosen Borel cross section $B_1\to X_1$ of the principle $L_1$-bundle $\pi_1:X_1\to B_1$
defines a Borel cocycle $\GA\times B_1\to L_1$, $(\ga,b)\mapsto\la_{\ga,b}$, by
\[
	\ga \la x_b=\la_{\ga,b}\la x_{\ga b}.
\]
The $m_1$-a.e. identity $p(\ga x)=\ga^\ta _*p(x)$ implies that  
for $\nu_1$-a.e. $b\in B_1$ and $m_{L_1}$-a.e. $\la\in L_1$
\[
	\ga^\ta_*p(\la x_b)=p(\ga \la x_b).
\]
Thus for $\ga\in\GA$ and $\nu_1$-a.e. $b\in B_1$ we have
\bqrn
	\ga^\ta_*\ps(b)&=&\lm \int_{L_1} \ga^\ta_*p(\la x_b)\,d\si_n(\la)
		=\lm \int_{L_1} p(\ga \la x_b)\,d\si_n(\la)\\
		&=&\lm \int_{L_1} p(\la_{\ga,b}\la x_{\ga b})\,d\si_n(\la)\\
		&=&\lm \int_{L_1} p(\la x_{\ga b})\,d\si_n(\la)=\ps(\ga b)
\eqrn
where the limits are understood with respect to the weak-* topology on
$\sP(B_2)$,  on which $\GA$ acts continuously.
The claim is proved.
\begin{clm}
There exists a measurable map $\ph:B_1\to B_2$ such that
$\nu_1$-a.e. $\ph(\ga b)=\ga^\ta \ph(b)$. 
Such a map is uniquely defined up to $\nu_1$-null sets and,
moreover, any $\ta$-twisted $\GA$-equivariant measurable 
map $\ps:(B_1,\nu_1)\to \sP(B_2)$ is $\nu_1$-a.e.
\[
	\ps(b)=\de_{\ph(b)}.
\]
\end{clm}
We shall use an existence of a $\ta$-twisted $\GA$-equivariant map
$\ps:(B_1,\nu_1)\to \sP(B_2)$ to construct a $\ta$-twisted $\GA$-equivariant 
map $\ph:(B_1,\nu_1)\to B_2$ (both defined up to $\nu_{1}$-null sets). 
Its uniqueness will follow from the construction (by alignment property)
and thereby will establish the essential uniqueness of $\ps(b)=\de_{\ph(b)}$.
 
Similarly to the map $p:X_1\to\sP(B_2)$ we define 
$q:X_2\to\sP(B_1)$ by $q_y={\pi_1}_*\nu_y$. 
By Lemma~\ref{L:pushforward} this is a measurable map satisfying 
\[
	q_{\ga^\ta y}=\ga_* q_y
\]
for all $\ga\in\GA$ and $m_2$-a.e. $y\in Y$.
Now let us use $\ps:B_1\to\sP(B_2)$ to pushforward $q_y$ and obtain $\et_y\in\sP(B_2)$
\[
	\et_y=\int_{B_1} \ps(b)\,dq_y(b),\qquad\quad
     \begin{CD}
      \sP(X_1)  @<{\nu_\cdot}<<  X_2 \\
       @VV{\pi_1}_*V    @VV{\et}V\\
        \sP(B_1) @>{\ps_*}>>  \sP(B_2)
     \end{CD}
\]
We need to invoke Lemma~\ref{L:pushforward} again to justify the following $m_2$-a.e. identity
\[
	\et_{\ga^\ta y}=\ps_*(q_{\ga^\ta y})=\ps_*(\ga_* q_y)=(\ga^\ta_*\ps)_*q_y=\ga^\ta_*\et_{y}.
\]
Indeed, if $E\in \sP(B_1)$ has $\nu_1(E)=0$ then $m_1(F)=0$ where $F=\pi_1^{-1}(E)$,
and $\nu_y(F)=0$ for $m_2$-a.e. $y\in X_2$, implying $q_y(E)=0$ for $m_2$-a.e. $y\in X_2$.

By the alignment property of $\pi_2:X_2\to B_2$ for $m_2$-a.e. $y\in X_2$
\[
	\de_{\pi_2(y)}=\et_y=\int_{B_1} \ps(b)\,dq_y(b)
\]
Dirac measures are extremal points of the convex compact set $\sP(B_2)$.
Hence we conclude that for $m_2$-a.e. $y\in X_2$
\begin{equation} \label{e:ps-qy}
	\ps(b)=\de_{\pi_2(y)}\qquad{\rm  for}\qquad q_y{\rm -a.e.}\  b\in B_1
\end{equation}
Thus $\ps(b)$ is a Dirac measure $\de_{\ph(b)}$ for $\nu_1$-a.e. $b\in B_1$. 
Previously established $\nu_1$-a.e. identity $\ps(\ga b)=\ga^\ta_* \ps(b)$ 
now implies 
\[
	\ph(\ga b)=\ga^\ta \ph(b)
\]
for $\nu_1$-a.e. $b\in B_1$ and all $\ga\in\GA$.
\begin{clm}
$\ph_*\nu_1\sim \nu_2$.
\end{clm}
Given a Borel subset $E _2\subset B_2$ let $E_1=\ph^{-1}(E_2)\subset B_1$.
Recall that by (\ref{e:ps-qy}) for $m_2$-a.e. $y\in X_2$ 
\begin{equation}\label{e:nuy}
	\nu_y\left(\Set{x\in X_1}{\ph(\pi_1(x))=\pi_2(y)}\right)=1.
\end{equation}
We have the following chain of equivalent conditions 
\bqrn 
    \nu_2(E_2)=0 &\iff& m_2(\Set{y\in X_2}{\pi_2(y)\in E_2})=0 \\
     & \iff & \nu_y(\Set{x\in X_1}{\pi_1(x)\in E_1})=0\quad m_2{\rm -a.e.}\ y \\
     & \iff & \bar{m}(\pi_1^{-1}E_1\times X_2)=0\\
     &\iff & m_1(\pi_1^{-1}E_1)=0\ \iff \ \nu_1(E_1)=0.
\eqrn
It proves that $\ph_*\nu_1\sim \nu_2$ as claimed.
\begin{clm}
$\ph:(B_1,[\nu_1])\to(B_2,[\nu_2])$ is a measurable measure class preserving bijection
which is $\nu_1$-a.e. $\ta$-twisted $\GA$-equivariant. 
The joining measure $\bar{m}$ is supported on $\sF_\ph$ as in 
Theorem~\ref{T:char-maps} (1).
\end{clm}
Since both alignment systems $(\pi_i:(X_i,m_i)\to B_i;\GA)$ are principle bundles
with amenable structure groups $L_i$ and compact $B_i$, by previous
considerations there exist ($\ta$ and $\ta^{-1}$-twisted) $\GA$-equivariant and 
measure class preserving maps
\[
	\ph:(B_1,[\nu_1])\to (B_2,[\nu_2]),\qquad 
	\ph^\prime: (B_2,[\nu_2])\to (B_1,[\nu_1]).
\]
The maps $\ps_1=\ph^\prime\circ\ph$ and $\ps_2=\ph\circ\ph^\prime$ are
well defined (up to $\nu_i$-null sets) measure class preserving self maps of $(B_i,\nu_i)$,
which are $\GA$-equivariant. 
By Lemma~\ref{L:align-uniq} it follows that $\ps_i(b)=b$ for $\nu_i$-a.e. $b\in B_i$.
We conclude that $\ph$ is a measurable measure class preserving bijection  
$(B_1,[\nu_1])\to (B_2,[\nu_2])$,
and along the way establish its uniqueness, mod null sets, as a $\ta$-twisted map.
Finally notice that (\ref{e:nuy}) actually means that $\bar{m}=\int \nu_y\,dm_2(y)$
is supported on the set $\sF_\ph$.

\begin{clm}
There exist compact subgroups $K_i< L_i$ so that the joining $m$ 
on $X_{1}\times X_{2}$ is $K_{1}\times K_{2}$-invariant; 
the quotient systems $(X_i^\prime,m_i^\prime)=(X_i,m_i)/K_i$ 
admit a $\ta$-twisted isomorphism $T^\prime:X_{1}^{\prime}\to X_{2}^{\prime}$ 
so that the original ergodic joining $\bar{m}$ descends to the measure
\[
	\int_{X_{1}^{\prime}}\delta_{x}\otimes\delta_{T^{\prime}x}\,dm_{1}^{\prime}(x)
\]
on the graph of $T^\prime$.
\end{clm}
Recall the disintegration (\ref{e:barm}) of the ergodic joining $\bar{m}$.
For $m_1$-a.e. $x\in X_1$ and $\mu_x$-a.e. $y\in X_2$ let 
$\et_{(x,y)}^{(2)}\in\sP(L_2)$ be the probability defined by 
\[
	\et_{(x,y)}^{(2)}(E)=\mu_x\Set{\la y}{\la\in E}\qquad E\in\sB(L_2).
\]
Denoting by $z=(x,y)$ points of $X_1\times X_2$ the above gives a measurable
map $z\mapsto \et_z^{(2)}\in\sP(L_2)$, well defined up to $\bar{m}$-null sets. 
Similarly there is a measurable map $z\mapsto \et_z^{(1)}\in\sP(L_1)$.

The $\GA$-action $\ga:(x,y)\mapsto (\ga x,\ga^\ta y)$ commutes with the $L_1\times L_2$-action. 
Hence for $i=1,2$ we have $\et^{(i)}_{\ga z}=\et^{(i)}_z$ for every $\ga\in\GA$ and $\bar{m}$-a.e. $z$. 
Hence $\et_z^{(i)}$ are $\bar{m}$-a.e. equal to fixed measures $\et_z^{(i)}=\et^{(i)}$ 
($i=1,2$). 

For both $i=1,2$ the probability measure $\et^{(i)}$ on $L_i$ are symmetric and
satisfy $\et^{(i)}=\et^{(i)}*\et^{(i)}$.
This follows from the fact that for $m_1$-a.e. $x\in X_1$ choosing $y$ and $y^\prime$ in $X_2$
independently according to $\mu_x$ we will have $y^\prime =\la y$ where the distribution of
$\la\in L_2$ is $\et^{(2)}$, and similarly for $\et^{(1)}$ (this is completely analogous to
the argument in the proof of part (3) of Theorem~\ref{T:generalCJQ}). 
Thus Lemma~\ref{L:idempotent-compact} yields that $\et^{(i)}$ is the normalized 
Haar measure on a \emph{compact subgroup} $K_i<L_i$.

Next consider the natural $\Gamma$-equivariant quotients 
\[
	p_{i}:(X_{i},m_{i})\to (X_i^\prime,m_i^\prime)=(X_i,m_i)/K_i\qquad(i=1,2)
\]
as in Lemma~\ref{L:quotient-by-comp}.
Consider $\sF_{\phi}$ with the joining measure $\bar{m}$ as an ergodic infinite measure preserving
action of $\Gamma$. The measure $\bar{m}$ is invariant under the action of the compact group 
$K=K_{1}\times K_{2}$ which commutes with the $\Gamma$-action. 
The quotient system $(\sF,\bar{m})/K$ is a subset of $X_{1}^{\prime}\times X_{2}^{\prime}$
with the measure $\bar{m}^{\prime}=\bar{m}/K$ having one-to-one projections on $X^{\prime}_{i}$.
Therefore it is supported on a graph of a measurable map $T^{\prime}:X_{1}^{\prime}\to X^{\prime}_{2}$,
and since $\bar{m}^{\prime}$ is invariant under the $\ta$-twisted diagonal $\Gamma$-action,
$T^{\prime}$ is a $\tau$-twisted isomorphism.

\bigskip

\noindent{Next assume that {\bf $K_{i}$ are normal in $L_{i}$} for both $i=1,2$.}

\medskip

Let $\Lambda_{i}=L_{i}/K_{i}$ and observe that $\pi_{i}:X_{i}\to B_{i}$ is a principle 
$\Lambda_{i}$-bundle which still has the alignemnt property (Lemma~\ref{L:align-interm}).
\begin{clm}
The groups $\Lambda_{i}=L_{i}/K_{i}$ are
continuously isomorphic.
\end{clm}
The graph of $T^\prime$ is supported on
\[
	\sF^\prime_\ph=\Set{(x_1^\prime,x_2^\prime)\in X^\prime_1\times 	X^\prime_2}{\ph(\pi_1^\prime(x_1^\prime))=\pi_2^\prime(x^\prime_2)}=\sF_{\phi}/(K_{1}\times K_{2}).
\]
This allows to define a Borel map $\rh:\LA_1\times X^\prime_1\to \LA_2$ by
\begin{equation}\label{e:Tprime-rho-cocycle}
	T^\prime(\la x)=\rh(\la,x) T^\prime(x)\qquad(\la\in \LA_1,\ x\in X_1^\prime)
\end{equation}
It is an (a.e.) cocycle i.e. $\rh(\la^\prime\la, x)=\rh(\la^\prime,\la x)\rh(\la,x)$ for $m_1^\prime$-a.e.
$x\in X^\prime_1$ and a.e. $\la\in \LA_1$.
Another a.e. identity is 
\[
	\rh(\la, \ga x)=\rh(\la,x)
\] 
for $\ga\in\GA$ whose action commutes 
with both $\LA_1$ on $X_1^\prime$ and $\LA_2$ on $X_2^\prime$.
Ergodicity implies that for a.e. $\la\in\LA_1$ the value $\rh(\la,x)$ is a.e. constant $\rh(\la)$.
The a.e. cocycle property of $\rh(\la,x)$ means that 
$\rh:\LA_1\to \LA_2$ is a measurable a.e. homomorphism.
It is well known that a.e. homomorphism between locally compact groups coincide a.e.
with a continuous homomorphism, which we continue to denote by $\rh$.
Thus (\ref{e:Tprime-rho-cocycle}) translates into (3) of Theorem~\ref{T:char-maps}, 
and the fact that $T^\prime_* m_1=m_2$ allows to conclude that $\rh:\LA_1\to\LA_2$ 
is one-to-one onto.

This completes the proof of Theorem~\ref{T:char-maps}.
\end{proof}

\section{Spaces of Horospheres}\label{S:Horospheres}
In this section we consider the geometric framework of pinched negatively curved manifolds.
Analysis of this geometric situation in particular implies the rigidity results for rank one symmetric 
spaces: Theorems~\ref{T:self-rk1} and \ref{T:rigidity-rk1} as we shall see in 
section~\ref{S:proofs-rigidity}.

\medskip

Let $N$ be a a complete simply connected Riemannian manifold with pinched negative 
curvature, let $\pa N$ denote its boundary at infinity (homeomorphic to a sphere $S^{\dim N-1}$), 
and $\GA<\Isom(N)$ be some non-elementary discrete group of isometries.

Recall some fundamental objects associated with this setup.
The critical exponent $\de=\de(\GA)$ of $\GA$ is 
\[
	\de(\GA)=\lim_{R\to\infty} \frac{1}{R} \log{\rm Card}\{\ga\in\GA \mid d(\ga p,p)<R\},
\]
(the limit exists and independent of $p$).
The Poincare series of $\GA$ is defined by
\[
	P_s(p,q)=\sum_{\ga\in\GA} e^{-s\cdot d(\ga p,q)}.
\] 
It converges for all $s>\de(\GA)$ and diverges for all $s<\de(\GA)$ regardless 
of the location of $p,q\in N$. 
If the series diverges at the critical exponent $s=\de(\GA)$ the group 
$\GA$ is said to be of \myemph{divergent type.}

Patterson-Sullivan measure/s is a measurable family $\{\nu_p\}_{p\in N}$ of 
mutually equivalent finite non-atomic measures supported on the limit set 
$L(\GA)\subs \pa N$ of $\Gamma$, satisfying
\begin{equation}\label{e:PS}
	\frac{d\nu_p}{d\nu_q}(\xi)=e^{-\de\cdot\be_\xi(p,q)}\rmand \nu_{\ga p}=\ga_*\nu_p
\end{equation}
where $\be_\xi(p,q)$ is the Busemann cocycle $\be_\xi(p,q)=\lim_{z\to\xi} \left[d(p,z)-d(q,z)\right]$
(the limit exists and is well defined for any $p,q\in N$ and $\xi\in\pa N$). 
Patterson-Sullivan measures exist, and for $\GA$ of divergent type the family 
$\{\nu_{p}\}$ is defined by the above properties uniquely, up to scalar multiple.

The space $\Hor(N)$ of horospheres is a principle $\bR$-bundle over 
$\pa N$ (see example~\ref{E:Hor}).
In the parametrization $\Hor(N)\iso \pa N\times \bR$ defined by a base point $o\in\ N$ 
the $\GA$-action is given by
\[
	\ga: (\xi,t)\mapsto (\ga\xi,t+c(\ga,\xi))\rmwhere c(\ga,\xi)=\be_{\xi}(\ga o,o).
\]
Define an infinite measure $m$ on $\Hor(N)$ by
\begin{equation}\label{e:PS-R}
	dm(\xi,t)=e^{-\de\cdot t}\,d\nu_{o}(\xi)\,dt
\end{equation}
where $\nu_{o}$ is the Patterson-Sullivan measure. 
Observe that $m$ is $\GA$-invariant.
We shall say that $\Gamma$ satisfies condition $({\bf E1})$ 

\begin{rem}\label{R:mp-extension}
Measure-theoretically $\Gamma$-action on $(\Hor(N), m)$ can be viewed as a standard
measure preserving extension of the measure class preserving $\Gamma$-action
on $(\pa N, \nu_{*})$ where $\nu_{*}$ is any representative of the measure class
of the Patterson-Sullivan measures.  
\end{rem}
%

\begin{thm}[Rigidity for actions of spaces of horospheres]
\label{T:rigidity-hor}
Let $N_1$ and $N_2$ be complete simply connected Riemannian manifolds of
pinched negative curvature, $\Isom (N_{1})>\Gamma_{1}\overto{\iso}\Gamma_{2}<\Isom(N_2)$ 
be two abstractly isomorphic discrete non-elementary torsion free groups of isometries 
with $\tau:\Gamma_{1}\overto{\iso}\Gamma_{2}$ denoting the isomorphism.  
Fis some base points $o_{i}\in N_{i}$ and let 
$\nu_{i}$ and $m_{i}$ denote the Patterson-Sullivan measures on 
$L(\Gamma_{i})\subseteq \pa N_{i}$ and their extensions on $\Hor(N_{i})$ 
respectively.

Assume that $(\Hor(N_{i},m_{i},\Gamma_{i})$ are ergodic. 
Then the following are equivalent:
\begin{enumerate}
\item
	The actions $(\Hor(N_{i}),m_{i},\GA_{i})$, $i=1,2$, admit an ergodic 
	${\rm Id}\times\ta$-twisted joining. 
\item
	The actions $(\Hor(N_{i}),m_{i},\GA_{i})$, $i=1,2$, admit a measurable
	$\ta$-twisted  isomorphism 
	\[
		T:(\Hor(N_{1}),m_{1})\overto{}(\Hor(N_{2}),m_{2}).
	\]
\item 
	There exists a ($\tau$-twisted) measure class preserving isomorphism 
	\[
		\phi: (L(\Gamma_{1}),\nu_{1})\overto{} (L(\Gamma_{2}),\nu_{2}).
	\]
\end{enumerate}
Under these (equivalent) conditions all ergodic joining are graphs of 
isomorphisms, the map $\phi$ is uniquely defined (up to null sets), and
every measurable isomorphism $\Hor(N_{1})\to \Hor(N_{2})$ has the form 
\[
	\hor_{\xi}(t)\in \Hor(N_{1})\quad \mapsto\quad \hor_{\phi(\xi)}(c\cdot t+s_{\xi}),
	\qquad \text{where}\qquad
	c=\delta_{2}/\delta_{1}.
\]

Furthermore, if $\Gamma_{i}$ are of divergent type, then the above conditions imply that 
$\phi:L(\Gamma_{1})\to L(\Gamma_{2})$ is a homeomorphism
and the groups $\GA_{i}$ have proportional length spectra
\[
	\de_{2}\cdot \ell_{2}(\ga^{\ta})=\de_{1}\cdot \ell_{1}(\ga)\qquad(\ga\in\GA_{1}),
\]
where $\ell_{i}(\ga)=\inf\{d_{i}(\ga p,p)\mid p\in N_{i} \}$ is the translation length 
of $\ga\in\Isom(N_{i})$.
\end{thm}
\begin{rem}\label{R:mp-ext}
In some cases the last statement implies that $N_{1}$ is equivariantly 
\emph{isometric} to $N_{2}$, after a rescaling by $\de_{1}/\de_{2}$. 
This is known as the Marked Length Spectrum Rigidity (Conjecture).
It has been proved for example for surface groups (\cite{Otal:MLSR:90}), 
or if say $N_{1}$ is  a symmetric space and $N_{1}/\Gamma_{1}$ is compact 
(\cite{Hamenstadt:cocycles:99}).
\end{rem}

\begin{proof}[Proof of Theorem~\ref{T:rigidity-hor}]

``(1) $\Rightarrow$ (2) and (3)''.
Recall that by Theorem~\ref{T:Hor} the $(\Hor(N_{i}),m_{i},\GA_{i})$
are alignment systems. 
Moreover these are $\bR$-principle bundles
with compact base spaces $\pa N_{i}$.
Let  $\bar{m}$ be an ergodic ${\rm Id}\times \tau$-twisted joining on 
$\Hor(N_{1})\times\Hor(N_{2})$.
Applying Theorem~\ref{T:char-maps} we get the desired measure class preserving
($\tau$-twisted) equivariant map
\[
	\phi:(\pa N_{1},\nu_{1})\overto{}(\pa N_{2},\nu_{2}).
\]
Its uniqueness being established on the way.
As $(\bR,+)$ has no compact subgroups we also conclude that $\bar{m}$
is a graph of an isomorphism. As the only automorphisms of $\bR$ are $t\mapsto c\cdot t$,
defining $s_{\xi}$ by $\hor_{\xi}(0)\mapsto \hor_{\phi(\xi)}(s_{\xi})$ we get the general
form of such an isomorphism as stated.

``(1) $\Rightarrow$ (2)'' and ``(3) or (2) $\Rightarrow$ (1)'' being trivial, we are left with 
``(3) $\Rightarrow$ (2)''.
But this follows from \ref{R:mp-ext}. 

\medskip

We are now left with the proof of the geometric conclusions in the case of divergent 
groups. The arguments below are probably known to experts, but we could not
find a good reference in the existing literature.  

We first recall some general facts from Patterson-Sulivan theory.
Let $\Gamma<\Isom(N)$ be a discrete group of isometries of a connected, 
simply connected manifold $N$ of pinched negative curvature.
Let 
\[
	\pa^2 N=\Set{(\xi,\et)\in\pa N\times\pa N}{\xi\neq\et}
\]
denote the space of pairs of distinct points at infinity of $N$ 
(this is the space of oriented but unparametrized geodesic lines in $N$).
Another Busemann cocycle (or Gromov product) can be defined for $\xi\neq\et\in\pa N$ 
and $p\in N$ by
\[
	B_p(\xi,\et)=\lim_{x\to\xi, y\to\et}\frac{1}{2}\left[d(p,x)+d(p,y)-d(x,y)\right].
\]
It can also be written as $B_p(\xi,\et)=\be_\xi(p,q)+\be_\et(p,q)$, 
where $q$ is an arbitrary point on the geodesic line $(\xi,\et)\subset N$. 
For any fixed $p\in N$ the function $B_{p}(\xi,\et)$ is continuous and proper on $\pa^{2} N$,
i.e. tends to $\infty$ as $\xi$ and $\et$ approach each other.

We have $\be_{\ga\xi}(\ga p,\ga q)=\be_{\xi}(p,q)$ and $B_{\ga p}(\ga\xi,\ga\et)=B_{p}(\xi,\et)$
for any isometry $\ga\in\Isom(N)$.
This implies that
\begin{equation}\label{e:B-be}
	B_p(\ga \xi,\ga \et)-B_p(\xi,\et)=\frac{1}{2}\left[\be_\xi(\ga p,p)+\be_\et(\ga p,p)\right].
\end{equation}
In view of (\ref{e:PS}) this means that the measure $\ol{\mu}$ on $\pa^{2} N\subset\pa N\times\pa N$ defined by
\begin{equation}\label{e:BPSPS}
	d\ol{\mu}(\xi,\et)=e^{2\de B_p(\xi,\et)}d\nu_p(\xi)\,d\nu_p(\et).
\end{equation}
is $\GA$-invariant.
This definition is independent of $p\in N$. 
One of the basic facts in Patterson-Sullivan theory states that $\GA$ is of divergent type iff 
its action on $(\pa^2 N,\ol{\mu})$ is ergodic (\cite{Sullivan:BAMS:82}, \cite{Yue:TAMS:96}).

Function $B_{o}(\cdot,\cdot)$ can also be used to define a \myemph{cross-ratio} on $\pa N$ by
\begin{equation}\label{e:cr-def}
	[\xi_{1},\xi_{2},\et_{1},\et_{2}]=e^{2\de\cdot\left[B_{o}(\xi_{1},\et_{1})+B_{o}(\xi_{2},\et_{2})
	-B_{o}(\xi_{1},\et_{2})-B_{o}(\xi_{2},\et_{1})\right]}
\end{equation}
where $o\in N$ is some reference point.
This cross-ratio is independent of the choice of $o\in N$, and  
is invariant under $\Isom(N)$ and in particular under $\GA$,
and satisfies the usual identities. 

\medskip

Returning to the given pair $\Gamma_{i}<\Isom(N_{i})$ we have
\begin{clm}
The measurable ($\tau$-twisted) equivariant map 
\[
	\ph:(L(\GA_{1},\nu_{1})\to (L(\GA_{2}),\nu_{2})
\] 
is a homeomorphism (possibly after an adjustment on null sets) and 
\[
	\ph(\ga \xi)=\ga^{\ta}\ph_{0}(\xi)
\] 
holds for all $\ga\in\GA_{1}$ and all $\xi\in L(\GA_{1})$.
Moreover
\begin{equation}\label{e:cr-ph}
	[\ph(\xi_{1}),\ph(\xi_{2}),\ph(\et_{1}),\ph(\et_{2})]_{2}
	=[\xi_{1},\xi_{2},\et_{1},\et_{2}]_{1}
\end{equation}
for all distinct $\xi_{1},\xi_{2},\et_{1},\et_{2}\in L(\GA_{1})\subs\pa N_{1}$.
\end{clm}
This is a consequence of property (i) and the ergodicity of $\GA_{i}$ on $\ol{\mu}_{i}$
(the following argument is a version of Sullivan's argument for Kleinian groups in 
\cite{Sullivan:BAMS:82}, 
but can also be traced back to Mostow in the context of quasi-Fuchsian groups).

The idea is that $\ph_*\nu_1\sim \nu_2$ and $\ol{\mu}_{i}\sim \nu_i\otimes\nu_i$ imply that the 
the pushforward measure $(\ph\times\ph)_* \ol{\mu}_1$ is absolutely continuous with 
respect to $\ol{\mu}_2$. 
Since $\ol{\mu}_{1}$ is $\GA_{1}$-invariant while 
$\ph$ is equivariant it follows that $(\ph\times\ph)_* \ol{\mu}_1$ is $\GA_{2}$-invariant.
Hence its Radon-Nikodym derivative with respect to $\ol{\mu}_{2}$ is a.e. a constant.
In view of (\ref{e:PS}), (\ref{e:BPSPS}) this amounts to a $\ol{\mu}_{1}$-a.e. relation    
\bqrn
	2\de_{2} \cdot B_{2}(\ph(\xi),\ph(\et))&=&2\de_{1}\cdot B_{1}(\xi,\et)+f(\xi)+f(\et)+C,\\
	\qquad
	\text{where}\qquad && f(\xi)=2\de_{1}\cdot\log\frac{d\ph_{*}\nu_{1}}{d\nu_{2}}(\ph(\xi)). 
\eqrn
Substituting these into the definition of the cross-ratios one observes
that the $f$-terms and the constant $C$ cancel out. It follows that
relation (\ref{e:cr-ph}) holds $\nu_{1}$-almost everywhere. 

For any fixed distinct $\xi_{2}, \xi_{3},\xi_{4}$ we have 
\[
	[\xi,\xi_{2}; \xi_{3},\xi_{4}]_{1}\to 0\qquad{\rm iff}\qquad
	\xi\to \xi_{3}.
\]
This allows, using Fubini theorem and a.e. identity (\ref{e:cr-ph}), to conclude that
$\ph$ agrees $\nu_{1}$-a.e. with a \emph{continuous} function $\ph_{0}$ defined 
on $\supp(\nu_{1})=L(\GA_{1})$. 
Since all the data are symmetric, it follows that
$\ph_{0}$ is a homeomorphism, relation (\ref{e:cr-ph}) extends from a.e. to everywhere
on $\supp(\nu_{1})=L(\Gamma_{1})$ by continuity. 

It is well known that the cross-ratio as above determines the marked length spectrum.
More precisely:
\begin{lem}
Let $N$ be a simply connected Riemannian manifold of pinched negative curvature, 
$\Gamma<\Isom(N)$ a non-elementary discrete group of isometries, 
$\delta=\delta(\Gamma)$ is growth exponent, 
and $[,;,]$ the corresponding cross ratio as in (\ref{e:cr-def}).
If $\gamma\in\Gamma$ is a hyperbolic element with attracting, repelling points
$\gamma_{+},\gamma_{-}\in\pa N$, then:
\[
	2\delta\cdot \ell(\gamma)=\log [\gamma_{+},\gamma_{-}; \xi,\gamma \xi]
\] 
 for all $\xi\in\pa N\setminus \{\gamma_{-},\gamma_{+}\}$.
\end{lem}
\begin{proof}
It is well known that $B_{p}(\xi,\eta)$ is within a constant (depending only on $N$) from
$\dist(p,(\xi,\eta))=\inf \{ d(p,x) \mid x\in (\xi,\eta)\}$.
For a fixed $\xi\neq \gamma_{\pm}$ and $p\in N$ we can estimate 
(with an error depending on $p$ and $\xi$, but independent of $n\in{\bf N}$) 
\[
	\dist(p,(\gamma_{+},\gamma^{n}\xi))=\dist(\gamma^{-n}p,(\gamma_{+},\xi))\asymp
	d(\gamma^{-n}p,p)\asymp n\cdot\ell(\gamma).
\]
Hence
\[
	\frac{1}{n} B_{p}(\gamma_{+},\gamma^{n}\xi)\overto{} \ell(\gamma).
\]
At the same time $\dist(p,(\gamma_{-},\gamma^{n}\xi))\overto{} \dist(p,(\gamma_{-},\gamma_{+}))$,
and so
\[
	\frac{1}{n} B_{p}(\gamma_{-},\gamma^{n}\xi)\overto{} 0.
\]
Since $\gamma$ fixes the points $\gamma_{-},\gamma_{+}\in\pa N$ and preserves the cross-ration we have for each $n$:
\bqrn
	&&\log[\gamma_{+},\gamma_{-};\xi,\ga\xi]=\frac{1}{n}\sum_{k=0}^{n-1}
		\log [\gamma_{+},\gamma_{-}\,;\,\ga^{k}\xi,\ga^{k+1}\xi]\\
	&&\quad =\frac{2\de}{n}\cdot \sum_{k=0}^{n-1} \left(B_{o}(\gamma_{+},\gamma^{k}\xi)-B_{o}(\gamma_{+},\gamma^{k+1}\xi)
			+B_{o}(\gamma_{-},\gamma^{k+1}\xi)-B_{o}(\gamma_{-},\gamma^{k}\xi)\right)\\
	&&\quad =\frac{2\de}{n}\cdot\left(B_{o}(\gamma_{+},\xi)-B_{o}(\gamma_{-},\xi)
		+B_{o}(\gamma_{+},\ga^{n}\xi)
	-B_{o}(\gamma_{-},\ga^{n}\xi)\right)\\
	&&\quad\overto{} \quad\ell(\gamma)\qquad\text{as}\quad n\to \infty.
\eqrn
\end{proof}

We return to $\Gamma_{i}<\Isom_{+}(N_{i})$, $i=1,2$, related by an abstract isomorphism  
$\tau:\Gamma_{1}\to\Gamma_{2}$
and an equivariant homeomorphism $\phi:L(\Gamma_{1})\to L(\Gamma_{2})$. 
The classification of elements of $\Isom_{+}(N_{i})$ into elliptic, parabolic and hyperbolic can be done
in terms of the dynamics on the boundaries, e.g. a hyperbolic isometry $g$ has two fixed points
$g_{-}$, $g_{+}$ and source/sink dynamics. 
Thus the topological conjugacy $\phi$ of the $\Gamma_{i}$-actions on $L(\Gamma_{i})\subset\pa N_{i}$,
shows that $\tau$ preserves the types of the elements.
Thus if $\gamma\in\Gamma_{1}$ is hyperbolic then so is $\gamma^{\tau}\in\Gamma_{2}$, and 
$\phi$ maps the corresponding repelling contracting points $\gamma_{\pm}$ of $\gamma$ 
to those of $\gamma^{\tau}\in\Gamma_{2}$, because 
\[
	\phi(\gamma_{\pm})=\phi(\lim_{n\to\pm\infty} \gamma^{n}\xi)=
	\lim_{n\to\pm\infty}(\gamma^{\tau})^{n}\phi(\xi)=
	\gamma^{\tau}_{\pm}
\]
for any $\xi\in\pa N_{1}\setminus \{\gamma_{-},\gamma_{+}\}$.
Thus using the previous Lemma we get:
\bqrn
	\delta_{1}\ell_{1}(\gamma)&=&\log[\gamma_{+},\gamma_{-};\xi,\ga\xi]_{1}
	=\log[\phi(\gamma_{+}),\phi(\gamma_{-});\phi(\xi),\phi(\ga\xi)]_{2}\\
	&=&\log[\gamma^{\tau}_{+},\gamma^{\tau}_{-};\phi(\xi),\ga^{\tau}\phi(\xi)]_{2}
	=\delta_{2}\ell_{2}(\gamma^{\tau}).
\eqrn
Hence $\delta_{1}\ell_{1}(\gamma)=\delta_{2}\ell_{2}(\gamma^{\tau})$ for all hyperbolic elements $\gamma\in\Gamma$ and the same formula (in the trivial form of $0=0$) applies to parabolic and elliptic $\gamma\in \Gamma_{1}$.
This completes the proof of Theorem~\ref{T:rigidity-hor}. 
\end{proof}

\section{Proofs of the Rigidity Results} \label{S:proofs-rigidity}
In this section we put all the ingredients developed in 
Sections~\ref{S:generalCJQ}-\ref{S:Horospheres} 
together in order to deduce the results stated in the Introduction.
We shall need some auxiliary facts, some of which, e.g. \ref{T:GP}, 
maybe of independent interest.
We start with the following general 
\begin{lem}
\label{L:pinched}
Let $\Gamma$ be some discrete group with ${\rm II}_{\infty}$ actions on six infinite measure
spaces linked into two sequences as follows:  
\[
	(\check{X}_{i},\check{m}_{i})\overto{p_{i}}(X_{i},m_{i})\overto{q_{i}}(\hat{X}_{i},\hat{m}_{i})
	\qquad(i=1,2).
\] 
Suppose that $\bar{m}$ be an ergodic (possibly $\tau$-twisted) joining of the $\Gamma$-actions
on $(X_{1},m_{1})$ with $(X_{2},m_{2})$. Then there exist ergodic ($\tau$-twisted) joinings
$\ol{\check{m}}$ of $(\check{X}_{1},\check{m}_{1})$ with $(\check{X}_{2},\check{m}_{2})$,
and $\ol{\hat{m}}$ of $(\hat{X}_{1},\hat{m}_{1})$ with $(\hat{X}_{2},\hat{m}_{2})$ 
so that 
\begin{equation}
\label{e:pinched}
	(\check{X}_{1}\times \check{X}_{2},\ol{\check{m}})\ 
	\overto{p_{1}\times p_{2}}\ 
	(X_{1}\times X_{2},\ol{m})\ 
	\overto{q_{1}\times q_{2}}\ 
	(\hat{X}_{1}\times\hat{X}_{2},\ol{\hat{m}})
\end{equation}
are quotient maps for the  ${\rm II}_{\infty}$ diagonal ($\Id\times\tau$-twisted) $\Gamma$-actions.
\end{lem}
\begin{proof}
The measure $\ol{\hat{m}}$ of $\ol{m}$ is defined by 
$\ol{\hat{m}}(E_{1}\times E_{2})=\ol{m}(q_{1}^{-1}E_{1}\times q_{2}^{-1}E_{2})$
and it is straightforward to verify that it is a joining of $\hat{m}_{1}$ with $\hat{m}_{2}$;
its ergodicity follows from the ergodicity of $\ol{m}$.

To construct $\ol{\check{m}}$ first consider the disintegration of $\check{m}_{i}$ 
with respect to $m_{i}$: 
\[
	\check{m}_{i}=\int_{X_{i}}\mu_{x}^{(i)}\,dm_{i}(x)\qquad(i=1,2).
\]
Consider the measure $\ol{m}^{*}$ on $\check{X}_{1}\times\check{X}_{2}$ defined by
\[
	\ol{m}^{*}=\int_{X_{1}\times X_{2}}\mu^{(1)}_{x}\otimes\mu^{(2)}_{y}\,dm(x,y).
\]
This measure forms a ($\Id\times\tau$-twisted) joining of the $\Gamma$-actions
on $(\check{X}_{i},\check{m}_{i})$ for $i=1,2$, and also projects to $\ol{m}$ under
$p_{1}\times p_{2}$.
Let $\ol{m}^{*}=\int \ol{\check{m}}_{t}\,d\eta(t)$ denote the ergodic decomposition
of $\ol{m}^{*}$ into ergodic joinings. 
Then $\ol{m}$ is an average of ergodic joinings $(p_{1}\times p_{2})_{*}\ol{\check{m}}_{t}$.
Since $\ol{m}$ is ergodic, $\eta$-a.e. ergodic joining $\ol{\check{m}}_{t}$ projects
to a multiple of $\ol{m}$ and can serve as $\ol{\check{m}}$ in the Lemma.  
\end{proof}

\bigskip

\subsection*{Proof of Theorem~\ref{T:self-rk1}}

Let $G$ be a real, connected, simple, non compact, center free, rank one group,
$G=KP$ and $P=MAN$ the Iwasawa decompositions.
Denote by $\bH=G/K$ the associated symmetric space and by $\pa\bH=G/P=G/MAN$
its boundary.
The unit tangent bundle is $S\bH=G/M$, and 
the space of horospheres $\Hor(\bH)$ can be identified with $G/MN$.

Let $H< G$ be a closed, unimodular, proper subgroup containing $N$,
and denote $\hat{H}=MN$, $\check{X}=G/N$, $X=G/H$, $\hat{X}=G/\hat{H}=\Hor(\bH)$
and let $\check{m}$, $m$, $\hat{m}$ denote the corresponding Haar measures.
We assume that $\Gamma$ acts ergodically on $(\check{X},\check{m})=(G/N,m_{G/N})$
and hence on $(X,m)$.
The projection
\[
	X=G/H\overto{} G/P, \qquad gH\mapsto gP
\] 
has the alignment property by Theorem~\ref{T:Hor} and Lemma~\ref{L:align-comp-ext}.
If $H$ is normal in $P$, for example if $H=\check{H}=N$ or $H=\hat{H}=MN$ then
$(X,m)\to B$ is a principle bundle with an alignment property (Theorem~\ref{T:Hor}),
and therefore Theorem~\ref{T:self-rk1} is a direct corollary of Theorem~\ref{T:generalCJQ}.

In the general case, the argument for algebraicity of quotients is the simplest:
any quotient $q:(X,m)\to (Y,n)$ defines a quotient of 
$(\check{X},\check{m})\overto{p} (X,m)\overto{q}(Y,n)$.
As was mentioned above Theorem~\ref{T:generalCJQ} applies to 
$\check{X}=G/\check{H}$, which gives that $(Y,n)$ can be identified with
$(G/H',m_{G/H'})$ where $\check{H}<H'$ with $H'/\check{H}$ compact.
Since the factor map $q\circ p: g\check{H}\mapsto gH'$ factors through $G/H$
it follows that $H<H'$ and $p(gH)=gH'$.

To analyze the centralizers of the $\Gamma$-action on $(X,m)$ we first 
consider general ergodic self joinings $\ol{m}$ of $(X,m)$.
Let $\ol{\check{m}}$ and $\ol{\hat{m}}$ be ergodic
self joinings of $(\check{X},\check{m})$ and $(\hat{X},\hat{m})$
as provided by Lemma~\ref{L:pinched}.
Applying Theorem~\ref{T:generalCJQ} to $(\check{X},\check{m})$
and $(\hat{X},\hat{m})$
we deduce that there exist $\lambda\in \check{\Lambda}=N_{G}(N)/N=P/N=MA$ and 
$a\in \hat{\Lambda}=N_{G}(MN)=P/MN=A$ so that
\[
	\ol{\check{m}}=const\cdot \int_{\check{X}}\delta_{x}\otimes \delta_{\lambda x}\,d\check{m}(x)
	\qquad\text{and}\qquad 
	\ol{\hat{m}}=const\cdot \int_{\hat{X}}\delta_{x}\otimes \delta_{a x}\,d\hat{m}(x). 
\]
Since $(X\times X,\ol{m})$ is an intermediate quotient
as in (\ref{e:pinched}), it follows that $a\in A$ is the image of $\lambda\in MA$ under
the natural epimorphism $\check\Lambda=MA\to\hat\Lambda=A$.
This completes the description of self-joinings in the theorem.

Finally, let $T:X\to X$ be a measurable centralizer of the $\Gamma$-action on $(X,m)$.
Applying the above arguments to the corresponding self-joining 
$\ol{m}=\int_{X}\delta_{x}\otimes \delta_{T(x)}\,dm(x)$ we deduce, in particular, that 
$T$ is ``covered'' by an algebraic automorphisms $\check{T}$ of $\check{X}=G/N$,
i.e. for some $q\in P=MAN$ the map $\check{T}:gN\mapsto  gqN$.
The fact that the graph of $\check{T}$ covers that of $T$ means that 
for a.e. $gH\in X$, the map $\hat{T}$ takes the preimage $p^{-1}(gH)=\{ghN\in \hat{X} \mid h\in H\}$
to $p^{-1}(g'H)=\{g'hN\in \hat{X} \mid h\in H\}$, where $T(gH)=g'H$. 
This implies that $q\in N_{G}(H)$, and $T(gH)=gqH$ -- an algebraic centralizer
as in \ref{E:algebraicCJQ}. 
Theorem~\ref{T:self-rk1} is proved.

\bigskip

\subsection*{Proof of Theorem~\ref{T:rigidity-rk1}}

Now consider a discrete subgroup $\GA< G$ which satisfies property (E2) 
(defined before the statement of Theorem~\ref{T:rigidity-rk1}). 
Such groups $\GA$ have the full limit set $L(\GA)=\pa\bH$, the maximal
critical exponent $\de(\GA)=\de(\bH)$ at which its Poincare series diverges
(so they are of divergent type),
and the associated Patterson-Sullivan measures are in the Haar class on $\pa \bH$.
The $\GA$-invariant measure $m$ on $S\bH\iso G/H_{1}$ as in (\ref{e:PS-R}) 
is a scalar multiple of $m_{G/H^{\prime}}$.
It is well known that $G=\Isom_{+}(\bH)$ can be identified with the \emph{conformal group}
on the boundary, the latter can be defined using the \emph{the cross-ratio}  
\begin{equation}\label{e:Isom-CR}
	G=\Isom_{+}(\bH)\iso\Set{\ps\in\Homeo_{+}(\pa\bH)}{[,,,]\circ\ps=[,,,]}
\end{equation}
Consider the framework of Theorem~\ref{T:rigidity-rk1} in which two rank-one groups $G_{i}$,
$i=1,2$, as above contain abstractly isomorphic discrete subgroups $\GA_{i}< G_{i}$,
$\ta:\GA_{1}\overto{\iso}\GA_{2}$, and the homogeneous space $X_{i}=G_{i}/H_{i}$
admit a $\ta$-twisted joining $\ol{m}$ with respect to the $\GA_{i}$-actions.
Denote by $\bH_{i}$, $\pa \bH_{i}$, $[,,,]_{i}$ etc. 
the corresponding symmetric spaces,
their boundaries,  cross-ratios etc.
Let set $\check{H}_{i}=N_{i}<H_{i}<\hat{H}_{i}=M_{i}N_{i}< G_{i}$ and 
\[
	\check{X}_{i}=G_{i}/\check{H}_{i}\ \overto{}\ 
	X_{i}=G_{i}/H_{i}\ \overto{}\ 
	\hat{X}_{i}=G_{i}/\hat{H}_{i}
	\qquad(i=1,2).
\]
Let $\ol{\hat{m}}$ on $\hat{X}_{1}\times \hat{X}_{2}$ denote the quotient joining 
of $\ol{m}$ on $X_{1}\times X_{2}$ as in Lemma~\ref{L:pinched}.  
Note that $\hat{X}_{i}=\Hor(\bH_{i})$.
Applying Theorems~\ref{T:Hor} and \ref{T:rigidity-hor} we conclude that 
there exists a homeomorphism $\ph:\pa\bH_{1}\to \pa\bH_{2}$ such that
\begin{itemize}
\item[(i)] 
	$[,,,]_{2}\circ\ph=[,,,]_{1}$ 
\item[(ii)]
	$\ph(\ga \xi)=\ga^\ta \ph(\xi)$ for all $\xi\in\pa\bH_{1}$, $\ga\in\GA_{1}$
\end{itemize}
In view of (\ref{e:Isom-CR}) property (i) yields an isomorphism 
$G_{1}\overto{\iso} G_{2}$ 
for which $\ph$ serves as the \emph{boundary map}. 
It follows from (ii) that this isomorphism extends $\tau:\Gamma_{1}\to\Gamma_{2}$.
Thus the result essentially reduces to that of Theorem~\ref{T:self-rk1},
(see the Proof of Theorem~\ref{T:rigidity-higher} for full details).

\subsection*{Proof of Theorem~\ref{T:self-higher}}
The proof of Theorem~\ref{T:self-rk1} applies almost verbatim to that of that of
Theorem~\ref{T:self-higher} with the appeal to Theorem~\ref{T:Hor} replaced by
Theorem~\ref{T:GH}.

\bigskip

For the proof of Theorem~\ref{T:rigidity-higher} we need some preparations,
which are of independent interest.

\medskip

Let $(B,\nu)$ be a standard probability space, $\Gamma$ a group acting by measure class 
preserving transformations on $(B,\nu)$. Such an action is called ``Strongly Almost Transitive'' if 
\[
	{\rm (SAT)}\qquad	
	\forall A\subset B\quad\text{with}\quad \nu(A)>0,\qquad\exists \gamma_{n}\in\Gamma:
	\quad\nu(\gamma_{n}^{-1}A)\to 1.
\]   
\begin{lem}
\label{L:SAT}
Let $\Gamma$ be a group with a measure class preserving (SAT) action on a
standard probability space $(B,\nu)$, let $C$ be a standard Borel space with a measurable
$\Gamma$-action, and $\pi_{1}, \pi_{2}:B\to C$ be two measurable maps such that:
\[
	\pi_{i}(\gamma x)=\gamma \pi_{i}(x)\quad\text{for}\qquad\nu\text{-a.e.}\ x\in B \qquad
	(\gamma\in\Gamma).
\]
Then $\pi_{1}(x)=\pi_{2}(x)$ for $\nu$-a.e. $x\in B$, unless the measures 
$(\pi_{1})_{*}\nu$, $(\pi_{2})_{*}\nu$ are mutually singular.
\end{lem}
\begin{proof}
Suppose that $\nu\left(\{x\in B\mid \pi_{1}(x)\neq \pi_{2}(x)\}\right)>0$. In this case there exists
a measurable set $E\subset C$ so that the symmetric difference 
$\pi_{1}^{-1}(E)\triangle \pi_{2}^{-1}(E)$ has positive $\nu$-measure.
Upon possibly replacing $E$ by its compliment, we may assume that the set
$A=\pi_{1}^{-1}(E)\setminus \pi_{2}^{-1}(E)$ has $\nu(A)>0$. 
Set $F_{i}=\pi_{i}(A)\subset C$. Then $F_{1}$ and $F_{2}$ are disjoint.
By the (SAT) property there exists a sequence $\{\gamma_{n}\}$ in $\Gamma$. so that
\[
	\sum_{n=1}^{\infty} \nu(\gamma_{n}^{-1}(B\setminus A))<\infty.
\] 
Then for $\nu$-a.e. $x\in B$ we have $g_{n}x\in A$ for all $n\ge n_{0}(x)\in\bN$.
For $i=1,2$ let $C_{i}\subset C$ denote the set of points $y\in C$ for which
$\{n \mid g_{n}y \not\in F_{i}\}$ is finite, i.e. $C_{i}=\liminf g_{n}^{-1}F_{i}$.
Hence $C_{1}\cap C_{2}=\emptyset$ because $F_{1}\cap F_{2}=\emptyset$.
In view of $\nu$-a.e. equivariance of $\pi_{i}$ we get that the measure $\eta_{i}=(\pi_{i})_{*}\nu$
is supported on $C_{i}$.  Hence $\eta_{1}\perp \eta_{2}$.
\end{proof}
\begin{thm}
\label{T:GP}
Let $G$ be a semi-simple group, $\Gamma<G$ a Zariski dense subgroup
with full limit set, e.g. a lattice in $G$ (irreducible or not),
$P<G$ a minimal parabolic and $\nu$ -- a probability measure on $G/P$
 in the Haar measure class.
  \begin{enumerate}
 \item
 	Let $Q<G$ be some parabolic subgroup containing $P$, 
	and $\pi:G/P\to G/Q$ be a measurable map, s.t. 
	$\pi(\gamma x)=\gamma \pi(x)$ a.e. on $G/P$ for all $\gamma\in\Gamma$.
	Then $\nu$-a.e. $\pi(gP)=gQ$.
 \item
 	Let $Q_{1}, Q_{2}<G$ be two parabolic subgroups containing $P$,
	$\nu_{i}$ probability measures on $G/Q_{i}$ in the Haar measure class,
	and $\phi:G/Q_{1}\to G/Q_{2}$ be a measurable bijection with $\phi_{*}\nu_{1}\sim\nu_{2}$
	and s.t. $\phi(\gamma x)=\gamma \phi(x)$ a.e. on $G/Q_{1}$ for all $\gamma\in\Gamma$.
	Then $Q_{1}=Q_{2}$ and $\phi(x)=x$ a.e. on $G/Q$
\end{enumerate}
\end{thm}
In particular the second statement with parabolic subgroups $Q_{1}=Q_{2}$ immediately gives:
\begin{cor}
Let $\Gamma<G$ be as in Theorem \ref{T:GP}, e.g. a lattice; 
let $Q<G$ be a parabolic subgroup and $\nu$ be a probability measure on $G/Q$ 
in the Haar measure class.
Then the measurable centralizer of the $\Gamma$-action on $G/Q$ is trivial.
\end{cor}
\begin{rem}
If $G$ is of higher rank, and $\Gamma<G$ is an irreducible lattice, then the 
only measurable quotients of the $\Gamma$-action on $G/Q$ are algebraic,
i.e. are given by $G/Q\to G/Q^{\prime}$ with $Q<Q^{\prime}$ and are given by 
$gQ\mapsto gQ^{\prime}$. 
This is the content of Margulis' Factor Theorem (see \cite{Margulis:book:89}).
\end{rem}
\begin{proof}[Proof of Theorem~\ref{T:GP}]
(1). The natural projection $\pi_{0}:G/P\to G/Q$, $\pi_{0}(gP)=gQ$ is $G$-equivariant.
It is well known that an action of a Zariski dense subgroup $\Gamma$ on $(G/P,\nu)$ is (SAT). 
The argument is then completed by Lemma~\ref{L:SAT}.

(2). Denote by $\pi_{i}:G/P\to G/Q_{i}$ the natural projections $\pi_{i}(gP)=gQ_{i}$.
Lemma~\ref{L:SAT} shows that the maps $\phi\circ \pi_{1}$ and $\pi_{2}$ agree $\nu$-a.e.
on $G/P$. This in particular implies that for a.e. $gP$ and all $q\in Q_{1}$:
\[
	gQ_{2}=\pi_{2}(gP)=\phi(\pi_{1}(gP))=\phi(gQ_{1})=\phi(gqQ_{1})=\pi_{2}(gqP)=gqQ_{2}.
\]
This means that $Q_{1}<Q_{2}$ and that $\phi(gQ_{1})=gQ_{2}$ a.e.
The same reasoning applies to $\phi^{-1}\circ\pi_{2}$, and $\pi_{1}$ as maps $G/P\to G/Q_{1}$,
giving $Q_{2}<Q_{1}$ and $\phi$ being a.e. the identity.
\end{proof}

\bigskip

\subsection*{Proof of Theorem~\ref{T:rigidity-higher}} 

Let $G$ be a semi-simple group, and $H_{i}<G$, ($i=1,2$) be two suitable subgroups
as in Definition~\ref{D:suitableH}. We consider the action of a lattice $\Gamma<G$
on the two homogeneous spaces $X_{i}=G/H_{i}$ equipped with the Haar measures
$m_{i}=m_{G/H_{i}}$.
 
Each of the the groups $H_{i}$ ($i=1,2$) is pinched between $\check{H}_{i}\normal\hat{H}_{i}$
associated to some parabolic subgroup $Q_{i}<G$.
The corresponding homogeneous spaces are linked by the natural $G$-equivariant
projections
\[
	\check{X}_{i}=G/\check{H}_{i}\ \overto{p_{i}}\ X_{i}=G/H_{i}\ \overto{q_{i}}\ \hat{X}_{i}=G/\hat{H}_{i}.
\]

These homogeneous spaces naturally project to $B_{i}=G/Q_{i}$, which is a compact space
with a continuous action of $G$. 
The $G$-equivariant projections:
\[
	\check{\pi}_{i}:\check{X}_{i}\overto{} B_{i},\qquad {\pi}_{i}:X_{i}\overto{} B_{i},
	\qquad \hat{\pi}_{i}:\hat{X}_{i}\overto{} B_{i}\qquad(i=1,2)
\]
have the alignment property with respect to the corresponding Haar measures
and the $\Gamma$-action (Theorem~\ref{T:GH}).
We also note that $\check{\pi}_{i}:\check{X}_{i}\to B_{i}$ and $\hat{\pi}_{i}:\hat{X}_{i}\to B_{i}$
are priciple bundles with amenable structure groups $\check{L}_{i}=Q_{i}/\check{H}_{i}$
and  $\hat{L}_{i}=Q_{i}/\hat{H}_{i}$.

Let $\ol{m}$ be an ergodic joining of the $\Gamma$-actions on $(X_{i},m_{i})$, and let 
$\ol{\check{m}}$ and $\ol{\hat{m}}$ be ergodic joinings of $\check{X}_{1}\times\check{X}_{2}$
and $\hat{X}_{1}\times\hat{X}_{2}$ as in Lemma~\ref{L:pinched}.
Applying Theorem~\ref{T:char-maps} to $\ol{\check{m}}$ (or $\ol{\hat{m}}$) we deduce
that there exists a $\Gamma$-equivariant measurable bijection 
$\phi:B_{1}\to B_{2}$ mapping the Haar measure class $[\nu_{1}]$ on $B_{1}$ 
to the measure class $[\nu_{2}]$ on $B_{2}$. 
Then part (2) of Theorem \ref{T:GP} shows that under these circumstances $Q_{1}=Q_{2}$.
Hence simplifying the notations we have a single parabolic $Q<G$ and
\[
	B<\check{H}<H_{1},H_{2}< \hat{H}\normal Q,
	\qquad \text{and}\qquad \check{L}=Q/\check{H},\qquad\hat{L}=Q/\hat{H}
\] 
and $\ol{\check{m}}$ and $\ol{\hat{m}}$ are ergodic self joinings of the $\Gamma$-actions
on $\check{X}=G/\check{H}$ and $\hat{X}=G/\hat{H}$, which are principle bundles over $B$
with the alignment property with respect to the $\Gamma$-action.
By Theorem~\ref{T:generalCJQ} they have the form:
\begin{equation}
\label{e:qinQ}
	\ol{\check{m}}=const\cdot \int_{G/\check{H}}\delta_{g\check{H}}\otimes 
		\delta_{gq\check{H}}\,d\check{m}(g\check{H}),
	\quad
	\ol{\hat{m}}=const\cdot \int_{G/\hat{H}}\delta_{g\hat{H}}\otimes 
		\delta_{gq\hat{H}}\,d\check{m}(g\hat{H})
\end{equation}
for some fixed $q\in Q$.

To complete the proof of the theorem, it remains to consider in more detail the case of an
isomorphisms $T:X_{1}\to X_{2}$; we shall prove that in this case
the above $q\in Q$ conjugates $H_{1}$ to $H_{2}$. 
So let $\ol{m}$ be the joining coming from the graph of $T$,
and let $\ol{\check{m}}$ be the ergodic joining supported on the graph of the map
$\check{T}(g\check{H})=g q \check{H}$ where $q\check{H}=\lambda\in\check{L}$.
Then $\check{T}$ maps the preimage 
\[
	p_{1}^{-1}(\{gH_{1}\})=\{gh_{1} \check{H}\in\check{X} \mid h_{1}\in H_{1}\}
\]
of a typical point $gH_{1}\in X_{1}$ to the preimage 
\[
	p_{2}^{-1}(\{g'H_{2}\})=\{g'h_{2} \check{H}\in\check{X} \mid h_{2}\in H_{2}\}
\]
of the point $g'H_{2}=T(g H_{1})\in X_{2}$. 
With $q\in Q$ as in (\ref{e:qinQ}) we have $gH_{1}q =g' H_{2}=T(gH_{1})$. 
This implies that $H_{2}=q^{-1}H_{1}q$ and $T(gH_{1})=gq H_{2}$ a.e.
This completes the proof of Theorem~\ref{T:rigidity-higher}.


\end{document}